\documentclass[11pt]{article}
\usepackage{amsmath}
\usepackage{amssymb}
\usepackage{latexsym}

\newtheorem{assumption}{Assumption}
\def\qed{ \ \vrule width.2cm height.2cm depth0cm\smallskip}

\newenvironment{proof}{\noindent {\bf Proof.\/}}{$\qed$\vskip 0.1in}

\newcommand{\la}{\langle}
\newcommand{\ra}{\rangle}

\newcommand{\ba}{\begin{array}}
\newcommand{\ea}{\end{array}}
\newcommand{\be}{\begin{equation}}
\newcommand{\ee}{\end{equation}}
\newcommand{\bea}{\begin{eqnarray}}
\newcommand{\eea}{\end{eqnarray}}
\newcommand{\beaa}{\begin{eqnarray*}}
\newcommand{\eeaa}{\end{eqnarray*}}

\def\dbD{\mathbb{D}}
\def\dbE{\mathbb{E}}
\def\dbF{\mathbb{F}}

\def\dbL{\mathbb{L}}
\def\dbM{\mathbb{M}}

\def\dbP{\mathbb{P}}
\def\dbR{\mathbb{R}}
\def\dbS{\mathbb{S}}

\def\a{\alpha}

\def\g{\gamma}
\def\d{\delta}
\def\e{\varepsilon}

\def\l{\lambda}

\def\si{\sigma}

\def\f{\varphi}

\def\o{\omega}

\def\O{\Omega}

\def\cE{{\cal E}}
\def\cF{{\cal F}}

\def\cH{{\cal H}}

\def\cL{{\cal L}}
\def\cM{{\cal M}}

\def\cP{{\cal P}}

\def\ms{\medskip}

\def\q{\quad}

\def\pa{\partial}
\def\cd{\cdot}
\def\cds{\cdots}

\def\tr{\hbox{\rm tr}}
\def\qed{ \hfill \vrule width.25cm height.25cm depth0cm\smallskip}

\newcommand{\basa}{\begin{assumption}}
\newcommand{\easa}{\end{assumption}}

\newcommand{\bas}{\begin{assum}}
\newcommand{\eas}{\end{assum}}

\def\limsup{\mathop{\overline{\rm lim}}}

\def\pa{\partial}

\def\cd{\cdot}
\def\cds{\cdots}

\def\tr{\hbox{\rm tr$\,$}}

\def\dis{\displaystyle}

\def\cad{{c\`{a}dl\`{a}g}}

\def\1{\mathbf{1}}

\def\:{\!:\!}
\def\reff#1{{\rm(\ref{#1})}}

at 9pt
\begin{document}

\title{\textbf{A Complete Representation Theorem for $G$-martingales }}
\author{Shige Peng \thanks{School of Mathematics, Shandong University, peng@sdu.edu.cn,
Research partially supported by NSF of China (No. 10921101) and National Basic Research Program
of China (973 Program) (No.2007CB814900).} \and Yongsheng Song\thanks{%
Academy of Mathematics and Systems Science, CAS, Beijing, China,
yssong@amss.ac.cn. Research supported by by NCMIS; Youth Grant of
National Science Foundation (No. 11101406); Key Lab of Random
Complex Structures and Data Science, CAS (No. 2008DP173182).}
\and Jianfeng Zhang\thanks{%
University of Southern California, Department of Mathematics,
jianfenz@usc.edu. Research supported in part by NSF grant DMS 10-08873.} }
\maketitle

\begin{abstract}
In this paper we establish a complete representation theorem for $G$%
-martingales. Unlike the existing results in the literature,
we provide the existence and uniqueness of the second order term, which
corresponds to the second order derivative in Markovian case. The main
ingredient of the paper is a new norm for that second order term, which is
based on an operator introduced by Song \cite{Song-unique}.
\end{abstract}

\date{}

\newtheorem{thm}{Theorem}[section] \newtheorem{lem}[thm]{Lemma} %
\newtheorem{cor}[thm]{Corollary} \newtheorem{prop}[thm]{Proposition} %
\newtheorem{rem}[thm]{Remark} \newtheorem{eg}[thm]{Example} %
\newtheorem{defn}[thm]{Definition} \newtheorem{assum}[thm]{Assumption}

\renewcommand {\theequation}{\arabic{section}.\arabic{equation}}



\section{Introduction}

\label{sect-Introduction} \setcounter{equation}{0} The notion of $G$-expectation is a type of nonlinear expectation proposed by Peng \cite{Peng-G, Peng-book}. In Markovian case, it corresponds to a fully nonlinear PDE. We also refer to Cheridito,  Soner,  Touzi and Victoir \cite{CSTV} and  Soner, Touzi and Zhang \cite{STZ-Duality, STZ-2BSDE}  for the closely related theory of Second Order Backward SDEs.  The theory has received very strong attention in the literature in recent years, we refer to the survey paper \cite{Peng-ICM} and the references therein, as well as some more recent developments: \cite{DNS}, \cite{ETZ}, \cite{HJPS}, \cite{MPZ},  \cite{Nutz}, \cite{NV}, \cite{NZ},  \cite{Song-unique}, to mention a few.
Their typical applications include, among others, stochastic optimization with diffusion control and  economic/financial models
with volatility uncertainty (see, e.g. \cite{DM}, \cite{EJ}, \cite{NS}) and numerical methods for high dimensional
fully nonlinear PDEs (see e.g. \cite{FTW}, \cite{Tan}, \cite{GZZ}).

$G$-expectation is a typical nonlinear expectation. It can be regarded as a
nonlinear generalization of Wiener probability space $(\Omega ,\mathcal{F}, \dbP_0)
$ where $\Omega =C([0,\infty) ,\mathbb{R}^{d})$, $\mathcal{F}=\mathcal{B}%
(\Omega )$ and $\dbP_0$ is a Wiener probability measure defined on $(\Omega ,%
\mathcal{F})$. Recall that the Wiener measure is defined such that the
canonical process $B_{t}(\omega ):=\omega _{t}$, $t\geq 0$ is a continuous
process with stable and independent increments, namely $(B_{t})_{t\geq 0}$
is a Brownian motion. $G$-expectation $\mathbb{E}^{G}$ is a sublinear
expectation on the same canonical space $\Omega $, such that the same
canonical process $B$ is a $G$-Brownian motion, i.e., it is a continuous
process with stable and independent increments. One important feature of
this notion is its time consistency. To be precise, let $\xi $ be a random
variable and $Y_{t}:=\mathbb{E}_{t}^{G}[\xi ]$ denote the conditional $G$%
-expectation, then one has $\mathbb{E}_{s}^{G}[\xi ]=\mathbb{E}_{s}^{G}[%
\mathbb{E}_{t}^{G}(\xi )]$ for any $s<t$. For this reason, we call the
conditional $G$-expectation a $G$-martingale, or a martingale under $G$%
-expectation. It is well known that a martingale under Wiener measure can be
written as a stochastic integral against the Brownian motion. Then a very
natural and fundamental question in this nonlinear $G$-framework is:
\begin{equation}
\label{Goal}
\mbox{What is the structure of a $G$-martingale $Y$?}
\end{equation}

Peng \cite{Peng-G} has observed that, for $Z\in \mathcal{H}_{G}^{2}$ and $%
\eta \in \mathcal{M}_{G}^{1}$ (see \reff{HpG} and \reff{MpG} below), the following process $Y$ is always a $G$-martingale:
\begin{equation}
dY_{t}=Z_{t}dB_{t}-G(\eta _{t})dt+{\frac{1}{2}}\eta _{t}d\langle B\rangle
_{t}.  \label{mrt3}
\end{equation}%
Here $G$ is the deterministic function Peng \cite{Peng-G} used to define $G$%
-expectations and $\langle B\rangle $ is the quadratic variation of the $G$%
-Brownian motion $B$. We remark that, in a Markovian framework, we have $%
Y_{t}=u(t,B_{t})$, where $u$ is a smooth function satisfying 
the following fully nonlinear PDE:
\begin{equation}
\partial _{t}u+G(\partial _{xx}u)=0.
\end{equation}%
Then $Z_{t}=\partial _{x}u(t,B_{t})$ and $\eta _{t}=\partial _{xx}u(t,B_{t})$%
. In particular, if $\xi =g(B_{T})$, then by PDE arguments we see
immediately that $Y_{t}:=\mathbb{E}_{t}^{G}[\xi ]$ has a representation
(\ref{mrt3}). Peng was even able to prove this $(Z,\eta )$%
-representation holds if $\xi $ is in a dense subspace
$\mathcal{L}_{ip}$ of $\mathcal{L}_{G}^{p}$ (see \reff{LpG} below). But observing that
$\mathcal{L}_{ip}$ is not a complete space, a very interesting
question was then raised to give a complete $(Z,\eta
)$-representation theorem for $\mathbb{E}_{t}^{G}[\xi ]$.

The first partial answer was provided by Xu and Zhang \cite{XZ}: if $Y$ is a
symmetric $G$-martingale, that is, both $Y$ and $-Y$ are $G$-martingales,
then
\begin{equation}
dY_{t}=Z_{t}dB_{t}\quad \mbox{for some process}~Z.
\end{equation}%
However, symmetric $G$-martingales captures only the linear part in this
nonlinear framework, and it is essentially important to understand the
structure of nonsymmetric $G$-martingales.

By introducing a new norm $\Vert \cdot \Vert _{\mathbb{L}_{G}^{2}}$ (see
\reff{dbLpG} below), Soner, Touzi and Zhang
\cite{STZ-G} proved a more general representation theorem: for $\xi \in
\mathbb{L}_{G}^{2}$,
\begin{equation}
dY_{t}=Z_{t}dB_{t}-dK_{t},
\label{mrt2}\end{equation}%
where $K$ is an increasing process such that $-K$ is a
$G$-martingale. It
has been proved independently in \cite{STZ-G} and Song \cite{Song-G} that $%
\mathbb{L}_{G}^{p}\supset \bigcap_{q>p}\mathcal{L}_{G}^{q}$, where $\Vert
\cdot \Vert _{\mathcal{L}_{G}^{q}}$ is the norm introduced in \cite{Peng-G}. Moreover,
\cite{Song-G} extended the representation \textrm{(\ref{mrt2})} to the case $%
p>1$.

Now the questions is, when does the process $K$ in \textrm{(\ref{mrt2})}
have the structure: $dK_{t}=G(\eta _{t})dt-{\frac{1}{2}}\eta _{t}d\langle
B\rangle _{t}$? Several efforts have been made in this direction. Hu and
Peng \cite{HP2} and Pham and Zhang \cite{PZ} made some progresses on the
existence of $\eta $. However, there is no characterization of the process $%
\eta $, and in particular, they do not provide an appropriate norm for $\eta
$. On the other hand, Song \cite{Song-unique} proved the uniqueness of $\eta
$ in the space $\mathcal{M}_{G}^{1}$. A clever operator was introduced in
this work, which successfully isolates the term ${\frac{1}{2}}\eta
_{t}d\langle B\rangle _{t}$ from $dK_{t}$, and thus essentially captures the
uncertainty of the underlying distributions. This idea turns out to be the
building block of the present paper.

Our main contribution of this paper is to introduce a norm for the process $%
\eta $, based on the work \cite{Song-unique}. We shall prove the existence and
uniqueness of the component $\eta$, which provides an essentially complete answer to Peng's question \reff{Goal}. 
Moreover, we shall provide a priori norm estimates. In particular, given $%
\xi _{1}$ and $\xi _{2}$ in appropriate space, let $(Y^{i},Z^{i},\eta ^{i})$%
, $i=1,2$, be the corresponding terms, we shall estimate the norms of $%
Z^{1}-Z^{2}$ and $\eta ^{1}-\eta ^{2}$ in terms of that of $Y^{1}-Y^{2}$,
where the latter one is more tractable due to the representation formula $%
Y_{t}=\mathbb{E}_{t}^{G}[\xi ]$. Unlike \cite{Song-unique}, we prove the
estimates via PDE arguments.

The rest of the paper is organized as follows. In Section \ref%
{sect-preliminary} we introduce the $G$-martingales and the involved spaces.
In Section \ref{sect-norm} we propose the new norm for $\eta$ and provide
some estimates. Finally in Section \ref{sect-mrt} we establish the complete
representation theorem for $G$-martingales.

\section{Preliminaries}

\label{sect-preliminary} \setcounter{equation}{0} In this section we
introduce $G$-expectations and $G$-martingales. We shall focus on a simple
setting in which we will establish the martingale representation theorem.
However, these notions can be extended to much more general framework, as in
many publications in the literature.

We start with some notations in multiple dimensional setting. Fix a
dimension $d$. Let $\dbR^d$ and $\dbS^{d}$ denote the sets of
$d$-dimensional column vectors and $d\times d$-symmetric matrices,
respectively.    For $\si_1, \si_2\in \dbS^d$, $\si_1
\le \si_2$ (resp. $\si_1 < \si_2$) means that $\si_2-\si_1$ is
nonnegative (resp. positive) definite, and we denote by $[\si_1,
\si_2]$ the set of $\si \in \dbS^d$ satisfying $\si_1 \le \si\le
\si_2$. Throughout the paper, we use ${\bf 0}$ to denote the $d$-dimensional  zero
vector or zero matrix, and $I_d$ the $d\times d$ identity matrix. For $x, \tilde x\in \dbR^d$, $\g, \tilde\g \in
\dbS^d$, define
\bea
\label{tr}
 x \cdot \tilde x := x^T \tilde{x},\q |x|:= \sqrt{x \cdot x}, &\mbox{and}&
  \g : \tilde \g := \tr (\g \tilde \g),\q  |\g| := \sqrt{\g:\g},
  \eea
  where $x^T$ denotes the transpose  of $x$. One can easily check that
   \bea
   \label{trest}
   |\g : \tilde\g| \le |\g||\tilde\g|, &\mbox{and}& -\g\le \tilde\g \le \g ~\mbox{implies that}~|\tilde\g|\le |\g|.
   \eea

\subsection{Conditional $G$-expectations}

We fix a finite time interval $[0,T]$, and two constant matrices $\mathbf{0}%
< \underline\sigma< \overline\sigma$ in $\mathbb{S}^d$. Define
\begin{eqnarray}  \label{G}
G(\gamma) := {\frac{1}{2}} \sup_{\sigma\in[\underline\sigma, \overline\sigma]%
} (\sigma^2: \gamma),\quad \mbox{for all}~\gamma\in \mathbb{S}^d.
\end{eqnarray}
Let $\Omega := \big\{\omega \in C([0,T], \mathbb{R}^{d}): \omega_0 = \mathbf{%
0}\big\}$ be the canonical space, $B$ the canonical process, and $\mathbb{F}%
:= \mathbb{F}^B$ the filtration generated by $B$. For $\xi = \varphi(B_T)$,
where $\varphi: \mathbb{R}^d \to \mathbb{R}$ is a bounded and Lipschitz
continuous function, following Peng \cite{Peng-G} we define the conditional $%
G$-expectation $\mathbb{E}^G_t[\xi] := u(t, B_t)$ where $u$ is the (unique)
classical solution of the following PDE on $[0,T]$:
\begin{eqnarray}  \label{PDE}
\partial_t u + G(\partial_{xx} u) =0,\quad u(T,x) = \varphi(x).
\end{eqnarray}
Let $\mathcal{L}_{ip}$ denote the set of random variables $%
\xi=\varphi(B_{t_1}, \cdots, B_{t_n})$ for some $0\le t_1<\cdots<t_n\le T$
and some Lipschitz continuous function $\varphi$. One may define $\mathbb{E}%
^G_t[\xi]$ in the same spirit, by defining it backwardly over each interval $%
[t_i, t_{i+1}]$. In particular, when $t=0$ we define $\mathbb{E}^G[\xi] :=
\mathbb{E}^G_0[\xi]$.

For any $p\ge 1$, define
\begin{eqnarray}  \label{LpG}
\|\xi\|_{\mathcal{L}^p_G}^p := \mathbb{E}^G[|\xi|^p],\quad \xi\in \mathcal{L}%
_{ip}.
\end{eqnarray}
Clearly this defines a norm in $\mathcal{L}_{ip}$. Let $\mathcal{L}^p_G$
denote the closure of $\mathcal{L}_{ip}$ under the norm $\|\cdot\|_{\mathcal{%
L}^p_G}$, taking the quotient as in the standard literature (i.e. we
do not distinguish random variables  $\xi_1$ and $\xi_2$ if
$\|\xi_1-\xi_2\|_{\mathcal{L}^p_G}=0$). As a mapping on the space
$\mathcal{L}_{ip}$, the conditional $G$-expectation is continuous
w.r.t.  the norm $\|\cdot\|_{\mathcal{%
L}^1_G}$. So one can easily extend it to all
$\xi\in\mathcal{L}^1_G$.

We next provide a representation of conditional $G$-expectations by
using the quasi-sure stochastic analysis, initiated by Denis and
Martini
\cite{DM} for superhedging problem under volatility uncertainty. Let $%
\mathcal{A}$ denote the space of $\mathbb{F}$-progressively measurable
processes taking values in $[\underline\sigma, \overline\sigma]$. Denoting
by $\mathbb{P}_0$ the Wiener measure, we define
\begin{eqnarray}  \label{cP}
\mathcal{P} := \Big\{ \mathbb{P}^\sigma := \mathbb{P}_0 \circ
(X^\sigma)^{-1}: \sigma\in \mathcal{A}\Big\} &\mbox{where}& X^\sigma_t :=
\int_0^t \sigma_s dB_s,~~\mathbb{P}_0\mbox{-a.s.}
\end{eqnarray}
Then $B$ is a $\mathbb{P}$-martingale for each $\mathbb{P}\in\mathcal{P}$.
Following \cite{DM}, we say
\begin{eqnarray}  \label{qs}
\mbox{a property holds $\cP$-quasi surely, abbreviated as $\cP$-q.s., if it
holds $\dbP$-a.s. for all $\dbP\in\cP$.}
\end{eqnarray}
We note that $\|\xi\|_{\cL^1_G} =0$ if and only if $\xi = 0$, $\cP$-q.s. Throughout this paper, random variables are considered the same if they are equal $\cP$-q.s. Then elements in $\cL^1_G$ can be viewed as standard random variables, but in  $\cP$-q.s. sense.  In particular, for any $\xi\in \cL^1_G$,  conditional $G$-expectation $\dbE^G_t[\xi]$ is  defined $\cP$-q.s.

It was proved in Denis, Hu and Peng \cite{DHP} that:
\begin{eqnarray}  \label{EGcP}
\mathbb{E}^G[\xi] = \sup_{\mathbb{P}\in\mathcal{P}} \mathbb{E}^\mathbb{P}%
[\xi],\quad \xi\in \mathcal{L}^1_G.
\end{eqnarray}
This result was extended by Soner, Touzi and Zhang \cite{STZ-G} to
conditional $G$-expectations: for any $\xi\in \cL^1_G$, $t\in [0,T]$, and $\mathbb{P}\in\mathcal{P}$,
\begin{eqnarray}  \label{EGcPt}
\mathbb{E}^G_t[\xi] = {\mathop{\rm ess\;sup}_{\mathbb{P}^{\prime }\in%
\mathcal{P}(t,\mathbb{P})}}^{\!\!\mathbb{P}}~ \mathbb{E}^{\mathbb{P}^{\prime
}}_t[\xi], ~~\dbP\mbox{-a.s.}, &\mbox{where}& \mathcal{P}(t,\mathbb{P}%
) := \Big\{\mathbb{P}^{\prime }\in\mathcal{P}: \mathbb{P}^{\prime }= \mathbb{%
P} ~\mbox{on}~\mathcal{F}_t\Big\}.
\end{eqnarray}
We remark that Peng \cite{Peng-g} had similar ideas, in the contexts of
strong formulation.

We finally note that $\mathbb{E}^G_t$ is obviously a sublinear expectation (again, all the equalities and inequalities are viewed in $\cP$-q.s. sense): for any $\xi, \xi_1, \xi_2\in \mathcal{L}^1_G$,
\bea  \label{sublinear}
\left.\ba{lll}
\mathbb{E}^G_t[\xi] = \xi,\q \mbox{if $\xi$ is $\cF_t$-measurable};&&  \mathbb{E}^G_t[\l \xi] =\l \xi,\q \mbox{for all $\l\ge 0$};\\
\mathbb{E}^G_t[\xi_1] \le \mathbb{E}^G_t[\xi_2],\q \mbox{if $\xi_1\le \xi_2$};&& \mathbb{E}^G_t[\xi_1 + \xi_2] \le \mathbb{E}^G_t[\xi_1] + \mathbb{E}%
^G_t[\xi_2].
\ea\right.
\eea

\subsection{Stochastic integrals}

First notice that, there exists a unique ($\cP$-q.s.)
$\mathbb{S}^d$-valued process $\langle B\rangle$ such that $B_t
B_t^T-\langle B\rangle_t$ is a symmetric $G$-martingale. In fact,
under each $\mathbb{P}\in \mathcal{P}$, $\langle B\rangle$ is the
same as the quadratic variation of the $\mathbb{P}$-martingale $B$,
and consequently,
\begin{eqnarray}  \label{<B>}
\underline\sigma^2 \le {\frac{d}{dt}}\langle B\rangle_t \le
\overline\sigma^2,\quad \mathcal{P}\mbox{-q.s.}
\end{eqnarray}
Naturally we call $\langle B\rangle$ the quadratic variation of $B$.
Next, we call an $\dbF$-progressively measurable  process $Z$ with appropriate dimension is an elementary process if it takes the form $Z= \sum_{i=0}^{n-1} Z_{t_i} \1_{[t_i, t_{i+1})}$ for some $0=t_0<\cds<t_n\le T$ and each component of $Z_{t_i}$ is in $\cL_{ip}$. Let $\cH^0_G $  denote the space of $\dbR^d$-valued elementary processes. For any $p\ge 1$, define
\begin{eqnarray}  \label{HpG}
\|Z\|_{\mathcal{H}^p_G}^p &:=& \mathbb{E}^G\Big[ \Big(\int_0^T (Z_tZ_t^T) :
d\langle B\rangle_t )\Big)^{\frac{p}{2}}\Big],\quad Z\in\mathcal{H}^0_G;
\end{eqnarray}
and let $\mathcal{H}^p_G$ denote the closure of $\mathcal{H}^0_G$ under the
norm $\|\cdot\|_{\mathcal{H}^p_G}$.

Now for each $Z\in \mathcal{H}^0_G$, we define its stochastic integral:
\begin{eqnarray}  \label{ZdB}
\int_0^t Z_s \cdot dB_s := \sum_{i=0}^{n-1} Z_{t_i}\cdot [B_{t_{i+1}\wedge
t} - B_{t_i\wedge t}],
\end{eqnarray}
One can easily prove the Burkholder-Davis-Gundy Inequality (see, e.g.
Song \cite{Song-G} Proposition 4.3): for any $p> 0$,
there exist constants $0< c_p < C_p<\infty$ such that
\begin{eqnarray}  \label{BDG}
c_p \|Z\|_{\mathcal{H}^p_G}^p \le \mathbb{E}^G\Big[\sup_{0\le t\le T}
|\int_0^t Z_s\cdot dB_s|^p\Big] \le C_p \|Z\|_{\mathcal{H}^p_G}^p.
\end{eqnarray}
Then one can extend the stochastic integral to all $Z\in\mathcal{H}^p_G$.

\subsection{$G$-martingales}

One important feature of conditional $G$-expectations is the time
consistency, which can also be viewed as dynamic programming principle:
\begin{eqnarray}  \label{timeconsistent}
\mathbb{E}^G_s\Big[\mathbb{E}^G_t(\xi)\Big] = \mathbb{E}^G_s[\xi],\quad%
\mbox{for all}\quad \xi\in \mathcal{L}^1_G ~~\mbox{and}~~0\le s<t\le T.
\end{eqnarray}
We recall that
\begin{eqnarray}  \label{Gmg0}
\mbox{a process $Y$ is called a $G$-martingale if $\dbE^G_s [Y_t] = Y_s$ for
all $0\le s<t\le T$.}
\end{eqnarray}
Therefore, $Y$ is a $G$-martingale if and only if $Y_t = \mathbb{E}^G_t[\xi]$
for $\xi = Y_T$.

Let $X, Y$ be two $G$-martingales. In general neither $-X$ nor $X+Y$
is a $G$-martingale since the conditional $G$-expectation is only
sublinear. If $-X$ is also a $G$-martingale, then we call $X$ a
symmetric $G$-martingale, and in this case one can easily check that $X+Y$ is still a
$G$-martingale.

It is clear that $\int_0^t Z_s \cd dB_s$ is a symmetric $G$-martingale
for all $Z \in
\mathcal{H}^1_G$. In particular, the canonical process $B$ is a symmetric $G$%
-martingale and is called a $G$-Brownian motion. However,
$G$-martingales have a richer structure. Let $\mathcal{M}^0_G$ be
the space of $\dbS^d$-valued elementary processes. Define
\begin{eqnarray}  \label{MpG}
\|\eta\|_{\mathcal{M}^p_G}^p &:=& \mathbb{E}^G\Big[\Big(\int_0^T |\eta_t| dt\Big)^p\Big]%
,\quad \eta\in \mathcal{M}^0_G;
\end{eqnarray}
and let $\mathcal{M}^p_G$ denote the closure of $\mathcal{M}^0_G$
under the norm $\|\cdot\|_{\mathcal{M}^p_G}$. An interesting fact
observed by Peng
\cite{Peng-G} is that the following decreasing process is also a $G$%
-martingale:
\begin{eqnarray}  \label{-K}
-K_t:= {\frac{1}{2}} \int_0^t \eta_s : d\langle B\rangle_s - \int_0^t
G(\eta_s) ds,\quad \eta\in \mathcal{M}^1_G.
\end{eqnarray}
Consequently, the following process $Y$ is always a $G$-martingale:
\begin{eqnarray}  \label{Gmg}
Y_t = Y_0 + \int_0^t Z_s\cdot dB_s - \Big[\int_0^t G(\eta_s) ds - {\frac{1}{2%
}} \int_0^t \eta_s : d\langle B\rangle_s\Big],\quad Z\in \mathcal{H}^1_G,
~\eta\in \mathcal{M}^1_G.
\end{eqnarray}

On the other hand, for any $\xi\in\mathcal{L}_{ip}$,  by
Peng \cite{Peng-book}  there exist $Z\in \mathcal{H}^1_G$ and $\eta\in
\mathcal{M}^1_G$ such that $Y_t := \mathbb{E}^G_t[\xi]$ satisfies \textrm{(%
\ref{Gmg})}. In particular, when $\xi = \varphi(B_T)$, for the classical
solution $u$ of PDE \textrm{(\ref{PDE})}, we have:
\begin{eqnarray}  \label{Y=u}
Y_t = u(t, B_t),\quad Z_t = \partial_x u(t, B_t),\quad \eta_t =
\partial_{xx} u(t, B_t).
\end{eqnarray}
Our goal of this paper is to answer the following natural question proposed
by Peng \cite{Peng-book}:
\begin{eqnarray}  \label{mrt}
\mbox{For what $\xi$ do there exist unique $Z\in \cH^1_G$ and $\eta\in
\cM^1_G$ satisfying \reff{Gmg}?}
\end{eqnarray}

The problem was partially solved by Soner, Touzi and Zhang \cite{STZ-G},
which introduced the following norm:
\begin{eqnarray}  \label{dbLpG}
\|\xi\|_{\mathbb{L}^p_G}^p := \mathbb{E}^G\Big[\sup_{0\le t\le T} \big(%
\mathbb{E}^G_t[|\xi|]\big)^p\Big],\quad \xi\in \mathcal{L}_{ip}.
\end{eqnarray}
Let $\mathbb{L}^p_G$ denote the closure of $\mathcal{L}_{ip}$ under the norm
$\|\cdot\|_{\mathbb{L}^p_G}$. Then for any $\xi\in \mathbb{L}^2_G$, there
exist unique $Z\in \mathcal{H}^2_G$ and an increasing process $K$ with $%
K_0=0 $ such that
\begin{eqnarray}  \label{mrt4}
Y_t := \mathbb{E}^G_t[\xi] = Y_0 + \int_0^t Z_s\cdot dB_s - K_t &\mbox{and}&
\|Z\|_{\mathcal{H}^2_G} + \|K_T\|_{\mathcal{L}^2_G} \le C\|\xi\|_{\mathbb{L}%
^2_G}.
\end{eqnarray}
It was proved independently by \cite{STZ-G} and Song \cite{Song-G} that $%
\|\xi\|_{\mathbb{L}^p_G} \le C_{p,q}\|\xi\|_{\mathcal{L}^q_G}$ for any $1\le
p<q$. Moreover, the above representation was extended by \cite{Song-G} to
the case $p>1$.

\subsection{Summary of notations}
For readers' convenience, we collect here some notations used in the paper:

\ms

$\bullet$ The inner product $\cdot$, the  trace operator $:$, and the norms $|x|$, $|\g|$ are defined by \reff{tr}.

$\bullet$ The function $G$, $G^\a$ and $G_\e$ are defined by \reff{G}, \reff{Ga}, and \reff{Ge}, respectively.

$\bullet$ The class of probability measures $\cP$, the $G$-expectation $\dbE^G$, and  the conditional $G$-expectation $\dbE^G_t$ are defined by \reff{cP}, \reff{EGcP}, and \reff{EGcPt} respectively.

$\bullet$ The norms $\|\xi\|_{\cL^p_G}$ and $\|\xi\|_{\dbL^p_G}$  for $\xi$ are defined by \reff{LpG} and \reff{dbLpG}, respectively.

$\bullet$ The norms $\|Z\|_{\cH^p_G}$ for $Z$ and $\|\eta\|_{\cM^p_G}$ for $\eta$ are defined by \reff{HpG} and \reff{MpG}, respectively.

$\bullet$ The norm $\|Y\|_{\dbD^p_G}$ for {\cad} processes $Y$, see also \reff{dbLpG}, is defined by:
\bea
\label{YpG}
\|Y\|_{\dbD^p_G}^p := \dbE^G\Big[\sup_{0\le t\le T} |Y_t|^p\Big].
\eea

$\bullet$ The operator $\cE^{\a}_{t_1, t_2}$ is defined by \reff{ua}.

$\bullet$ The constants $c_0, C_0$ are defined by \reff{c0}.

$\bullet$ The function $\d_n$ is defined by \reff{dn}.

$\bullet$ The new norms $\|\eta\|_{\dbM_G}$ and $\|\eta\|_{\dbM^*_G}$ for $\eta$ are defined by \reff{dbMG} and \reff{dbM*G}, respectively.

$\bullet$ The space $\cM^1_{G_0}$ and class $\cP_0$ are defined by \reff{cM1G0} and \reff{cPe}, respectively.

$\bullet$ The new metric $d_{G,p}(\xi_1, \xi_2)$ for $\xi$ is defined by \reff{rhoG}, and $\dbL^{*p}_G$ is the corresponding closure space.

$\bullet$ For $0 \le s\le t\le T$, the shifted canonical process $B^s_t$ is defined by:
\bea
\label{Bst}
B^s_t := B_t - B_s.
\eea

\section{A new norm for $\protect\eta$}

\label{sect-norm} \setcounter{equation}{0}

Our main contribution of the paper is to introduce a norm for
$\eta$. For that purpose, we shall introduce two nonlinear
operators, one via PDE arguments and the other via probabilistic
arguments. The latter one is strongly motivated by the work Song
\cite{Song-unique}, and the connection between the two operators is
established in Lemma \ref{lem-dB} below.

\subsection{The nonlinear operator via PDE arguments}

We first introduce a new nonlinear operator $\mathcal{E}^{\alpha}$ on
Lipschitz continuous functions, with a parameter $\alpha\in \mathbb{S}^d$.
Define
\begin{eqnarray}  \label{Ga}
G^\alpha(\gamma)={\frac{1}{2}}[G(\gamma+2\alpha)+G(\gamma-2\alpha)],\quad
\gamma\in \mathbb{S}^d.
\end{eqnarray}
Given $0\le t_1 < t_2 \le T$ and a Lipschitz continuous function $\varphi$,
define $\mathcal{E}^{\alpha}_{t_1, t_2}(\varphi) := u^\alpha(t_1, \cdot)$,
where $u^\alpha$ is the unique viscosity solution of the following PDE on $[t_1, t_2]$:
\begin{eqnarray}  \label{ua}
\partial_t u^\alpha + G^\alpha(\partial_{xx} u^\alpha) =0,\quad
u^\alpha(t_2,x) = \varphi(x).
\end{eqnarray}
Clearly $G^\alpha$ is strictly increasing and convex in $\gamma$. In
particular, the above PDE is parabolic and is wellposed. We collect below
some obvious properties of $G^\alpha$ and $\mathcal{E}^{\alpha}$, whose
proofs are omitted.

\begin{lem}
\label{lem-Ga} For any $\a\in \dbS^d$,

(i) $\cE^{\a}$ satisfies the semigroup property:
\bea
\label{semigroup}
\cE^{\a}_{t_1, t_2}\big(\cE^{\a}_{t_2, t_3}(\f)\big) = \cE^{\a}_{t_1, t_3}(\f),&\mbox{for any}& 0\le t_1<t_2 <t_3\le T.
\eea

(ii) $G^{-\a}=G^\a \ge G = G^{{\bf 0}} $.


(iii) If $\f = c$ is a constant, then $\cE^{\a}_{t_1, t_2}(c) = c + G^\a({\bf 0}) (t_2-t_1)$.
\end{lem}

The next property will be crucial for our estimates. Let
\begin{eqnarray}  \label{c0}
c_0 := \mbox{the smallest eigenvalue of}~ {1\over 2}[\overline\si^2 - \underline \si^2],  &\mbox{and}& C_0 :=  {1\over 2}|\overline\si^2 - \underline \si^2|.
\eea
Then clearly $C_0 \ge c_0 >0$ and $\underline \si^2  + c_0 I_{d} \le  \overline \si^2 -  c_0 I_{d}$. Denote, for $\e \le c_0$,
\begin{eqnarray}  \label{Ge}
G_\varepsilon(\gamma) := {\frac{1}{2}} \sup_{\sigma\in[\underline\sigma_%
\varepsilon, \overline\sigma_\varepsilon]} (\sigma^2: \gamma),&\mbox{where}%
&\underline\sigma_\varepsilon^2:= \underline\sigma^2 + \varepsilon I_d ,~~
\overline\sigma_\varepsilon^2:=\overline\sigma^2-\varepsilon I_d.
\end{eqnarray}
\begin{lem}
\label{lem-Gae}
(i) For any $0<\e\le c_0$ and $\a, \g\in \dbS^d$,  it holds that
\bea
\label{Gaest}
G_\e(\g) +  \e |\a| \le G^\a(\g) \le  G(\g) + C_0|\a|.
\eea

(ii) Assume $\underline \f \le \f \le \overline \f$ are Lipschitz continuous functions, and $0\le t_1 < t_2 \le T$. Then
\beaa
\dbE^{G_\e}\Big[ \underline \f (x+B_{t_2}^{t_1})\Big] + \e |\a|(t_2-t_1)  \le \cE^\a_{t_1, t_2}(\f) (x) \le \dbE^{G}\Big[ \overline \f (x+B_{t_2}^{t_1})\Big] + C_0 |\a|(t_2-t_1).
\eeaa
\end{lem}
{\noindent \textit{Proof.\quad}}  (i) We first
prove the left inequality. Let $\alpha_1, \cdots, \alpha_d$ denote the
eigenvalues of $\alpha$, and $\hat \alpha$ the diagonal matrix with
components $\alpha_1,\cdots, \alpha_d$. Then $|\alpha| = (\alpha_1^2 +
\cdots+\alpha_d^2)^{\frac{1}{2}}$, and there exists an orthogonal matrix $Q$
such that $Q^T\alpha Q = \hat\alpha $. Let $\hat c_\varepsilon$ denote a
diagonal matrix whose diagonal components take values $\varepsilon$ or $%
-\varepsilon$. Now for any $\sigma_\varepsilon \in
[\underline\sigma_\varepsilon, \overline\sigma_\varepsilon]$, by \textrm{(%
\ref{Ge})}, we have
\begin{eqnarray*}
\sigma_\varepsilon^2 + Q\hat c_\varepsilon Q^T\in [\underline \sigma^2,
\overline \sigma^2] &\mbox{and}& \sigma_\varepsilon^2 - Q\hat c_\varepsilon
Q^T \in [\underline \sigma^2, \overline \sigma^2].
\end{eqnarray*}
Then
\begin{eqnarray*}
2G^\alpha(\gamma) &=& G(\gamma+2\alpha) + G(\gamma-2\alpha) \\
&\ge& {\frac{1}{2}}\Big[(\sigma_\varepsilon^2 + Q\hat c_\varepsilon
Q^T):(\gamma + 2\alpha) + (\sigma_\varepsilon^2 - Q\hat c_\varepsilon
Q^T):(\gamma - 2\alpha)\Big] \\
&=& \sigma_\varepsilon^2 :\gamma + 2(Q\hat c_\varepsilon Q^T) :\alpha =
\sigma_\varepsilon^2 :\gamma + 2 \hat c_\varepsilon : (Q^T \alpha Q ) =
\sigma_\varepsilon^2 :\gamma + 2\hat c_\varepsilon : \hat\alpha.
\end{eqnarray*}
By the arbitrariness of $\sigma_\varepsilon$ and $\hat c_\varepsilon$, we
get
\begin{eqnarray*}
G^\alpha(\gamma) \ge G_\varepsilon(\gamma) + \varepsilon\sum_{i=1}^d
|\alpha_i| \ge G_\varepsilon(\gamma) + \varepsilon|\alpha|.
\end{eqnarray*}

We now prove the right inequality of \textrm{(\ref{Gaest})}. For any $%
\sigma_1, \sigma_2 \in [\underline\sigma, \overline\sigma]$, we have
\begin{eqnarray*}
\sigma_1^2 : (\gamma + 2\alpha) + \sigma_2^2 : (\gamma - 2\alpha) =
(\sigma_1^2 + \sigma_2^2) : \gamma + 2(\sigma_1^2 - \sigma_2^2): \alpha.
\end{eqnarray*}
Note that
\begin{eqnarray*}
\underline\sigma^2 \le {\frac{1}{2}} (\sigma_1^2 + \sigma_2^2) \le
\overline\sigma^2,\quad -[\overline\sigma^2 -
\underline\sigma^2]\le \sigma_1^2 - \sigma_2^2\le \overline\sigma^2 -
\underline\sigma^2.
\end{eqnarray*}
Then, by \textrm{(\ref{trest})},
\begin{eqnarray*}
\sigma_1^2 : (\gamma + 2\alpha) + \sigma_2^2 : (\gamma - 2\alpha) \le 4
G(\gamma) + 4C_0 |\alpha|.
\end{eqnarray*}
Since $\sigma_1, \sigma_2$ are arbitrary, we prove the right inequality of
\textrm{(\ref{Gaest})}, and hence \textrm{(\ref{Gaest})}.

(ii) One can easily check that
\begin{eqnarray*}
&&\mathbb{E}^{G_\varepsilon}\Big[ \underline \varphi (x+B_{t_2}^{t_1})\Big] +
\varepsilon |\alpha|(t_2-t_1) = \underline v^\alpha(t_1, x), \\
&& \mathbb{E}^{G}\Big[ \overline \varphi (x+B_{t_2}^{t_1})\Big] + C_0
|\alpha|(t_2-t_1)=\overline v^\alpha(t_1,x),
\end{eqnarray*}
where $\underline v^\alpha, \overline v^\alpha$ are the unique viscosity
solution of the following PDEs on $[t_1, t_2]$:
\begin{eqnarray*}
&&\partial_t \underline v^\alpha + G_\varepsilon(\partial_{xx} \underline
v^\alpha) + \varepsilon |\alpha| =0,\quad \underline v^\alpha(t_2, x) =
\underline \varphi(x); \\
&&\partial_t \overline v^\alpha + G(\partial_{xx} \overline v^\alpha) +
C_0|\alpha| =0,\quad \overline v^\alpha(t_2, x) = \overline \varphi(x).
\end{eqnarray*}
Then the statement follows directly from \textrm{(\ref{Gaest})} and the
comparison principle of PDEs. \hfill \vrule width.25cm height.25cm
depth0cm\smallskip

\subsection{The nonlinear operator via probabilistic arguments}

For any $n\ge 1$, denote $t^n_i := {\frac{i}{n}}T$, $i=0,\cdots,n$,
and define
\begin{eqnarray}  \label{dn}
\delta_n(t)=\sum_{i=0}^{n-1}(-1)^i1_{[t^n_i, t^n_{i+1})}, \quad t\in[0,T].
\end{eqnarray}
This function was introduced in \cite{Song-unique} which plays a key
role for constructing a new norm for process $\eta$. According to \cite{Song-unique},
we have
\begin{lem}
\label{lem-detadt}
For any $\eta\in \cM^1_G$, it holds that
$
\lim_{n\to \infty} \dbE^G\Big[\int_0^T G(\eta_t) \d_n(t) dt\Big] = 0.
$
\end{lem}

The next lemma establishes the connection between $\delta_n$ and $(G^\alpha,  \cE^{\a})$.
\begin{lem}
\label{lem-dB}
Let $0\le s < t\le T$ and $\a\in\dbS^d$.

(i)  For any $\g\in\dbS^d$, we have
\bea
\label{dg}
\lim_{n\to\infty}  \dbE^G_s\Big[\int_s^t [\a \d_n(r) +{1\over 2} \g] : d\la B\ra_r\Big] = G^\a(\g)(t-s).
\eea

(ii) For any $x\in\dbR^{d}$ and any Lipschitz continuous function $\f$, we have
\bea
\label{df}
\lim_{n\to\infty}  \dbE^G_s\Big[\int_s^t  \d_n(r) \a: d\la B\ra_r + \f(x+ B_t^s)\Big] = \cE^{\a}_{s,t}(\f)(x).
\eea
\end{lem}
{\noindent \textit{Proof.\quad}} (i) Fix $n$ such that ${2T\over n}<t-s$. Note that
\begin{eqnarray*}
&&\mathbb{E}^G_{t^n_{2i}}\Big[\int_{t^n_{2i}}^{t^n_{2i+2}} [\alpha
\delta_n(r) +{\frac{1}{2}}\gamma]: d\langle B\rangle_r\Big] \\
&=& \mathbb{E}^G_{t^n_{2i}}\Big[({\frac{1}{2}}\gamma+\alpha):[\langle
B\rangle_{t^n_{2i+1}} - \langle B\rangle_{t^n_{2i}}]+({\frac{1}{2}}\gamma -
\alpha ):[\langle B\rangle_{t^n_{2i+2}} - \langle B\rangle_{t^n_{2i+1}}]\Big]
\\
&=& \mathbb{E}^G_{t^n_{2i}}\Big[({\frac{1}{2}}\gamma+\alpha):[\langle
B\rangle_{t^n_{2i+1}} - \langle B\rangle_{t^n_{2i}}] + \mathbb{E}%
^G_{t^n_{2i+1}}\big[({\frac{1}{2}}\gamma - \alpha ):[\langle
B\rangle_{t^n_{2i+2}} - \langle B\rangle_{t^n_{2i+1}}]\big]\Big] \\
&=& \mathbb{E}^G_{t^n_{2i}}\Big[({\frac{1}{2}}\gamma+\alpha):[\langle
B\rangle_{t^n_{2i+1}} - \langle B\rangle_{t^n_{2i}}] +G(\gamma - 2\alpha ){%
\frac{T}{n}} \Big] \\
&=& \mathbb{E}^G_{t^n_{2i}}\Big[({\frac{1}{2}}\gamma+\alpha):[\langle
B\rangle_{t^n_{2i+1}} - \langle B\rangle_{t^n_{2i}}]\Big]  +G(\gamma - 2\alpha ){%
\frac{T}{n}} \\
&=& G(\gamma+2\alpha) {\frac{T}{n}} +G(\gamma - 2\alpha ){\frac{T}{n}} =
G^\alpha(\gamma)(t^n_{2i+2}-t^n_{2i}).
\end{eqnarray*}
Similarly, for any $i<j$,
\begin{eqnarray*}
\mathbb{E}^G_{t^n_{2i}}\Big[\int_{t^n_{2i}}^{t^n_{2j}}[\alpha \delta_n(r) +{%
\frac{1}{2}}\gamma]:d\langle B\rangle_r\Big] =
G^\alpha(\gamma)(t^n_{2j}-t^n_{2i}).
\end{eqnarray*}
Now assume $t^n_{2i} \le s < t^n_{2i+1} \le t^n_{2j} \le t < t^n_{2j+2}$.
Then
\begin{eqnarray*}
&& \Big| \mathbb{E}^G_s\Big[\int_s^t [\alpha \delta_n(r) +{\frac{1}{2}}%
\gamma] : d\langle B\rangle_r\Big] - G^\alpha(\gamma) (t-s) \Big| \\
&\le&\Big|\mathbb{E}^G_s\Big[\int_s^t [\alpha \delta_n(r) +{\frac{1}{2}}%
\gamma] : d\langle B\rangle_r\Big] - \mathbb{E}^G_s\Big[%
\int_{t^n_{2i+2}}^{t^n_{2j}} [\alpha \delta_n(r) +{\frac{1}{2}}\gamma] :
d\langle B\rangle_r\Big] \Big| \\
&& + \Big|G^\alpha(\gamma)(t^n_{2j}-t^n_{2i+2})-G^\alpha(\gamma) (t-s) \Big|
\\
&\le& \mathbb{E}^G_s\Big[ \Big|\lbrack \int_s^{t^n_{2i+2}} +
\int_{t^n_{2j}}^t ][\alpha \delta_n(r) +{\frac{1}{2}}\gamma]:
d\langle
B\rangle_r\Big| \Big] + {\frac{2T}{n}}|G^\alpha(\gamma)| \\
&\le& {\frac{2T}{n}}|\overline\sigma^2| [|\alpha|+{\frac{1}{2}}|\gamma|] + {%
\frac{2T}{n}}|G^\alpha(\gamma)|\to 0,\quad\mbox{as}~n\to\infty,
\end{eqnarray*}
where the last inequality thanks to  \reff{trest}.
This proves the result.

(ii) Without loss of generality, assume $t=T$. Define
\begin{eqnarray*}
\overline{u}(t,x) &:=&\mathop{\overline{\rm lim}}_{n\rightarrow\infty}
\overline{u}^n(t,x) := \mathop{\overline{\rm lim}}_{n\rightarrow\infty}
\mathbb{E}^G_t\Big[\int_t^T \delta_n(r) \alpha: d\langle B\rangle_r +
\varphi(x+ B_T^t)\Big], \\
\underline{u}(t,x)&:=&\mathop{\overline{\rm lim}}_{n\rightarrow\infty}
\underline{u}^n(t,x) := \mathop{\underline{\rm lim}}_{n\rightarrow\infty}%
\mathbb{E}^G_t\Big[\int_t^T \delta_n(r)\alpha: d\langle B\rangle_r +
\varphi(x+ B_T^t)\Big].
\end{eqnarray*}
By the structure of $G$-framework it is clear that $\underline{u}$ and $%
\overline{u}$ are deterministic functions. Obviously $\underline{u}\le
\overline{u}$. We claim that $\overline{u}$ and $\underline{u}$ are
viscosity subsolution and viscosity supersolution of PDE \textrm{(\ref{ua})}
with $t_1=0, t_2=T$. Note that PDE \textrm{(\ref{ua})} satisfies the
comparison principle for viscosity solutions. Then $\overline{u}\le
\underline{u}$ and thus $\overline{u}(t,x)= \underline{u}(t,x) = \mathcal{E}%
^{\alpha}_{t,T}(\varphi)(x) $. This proves the result.

We now prove that $\overline{u}$ is a viscosity subsolution, and the
viscosity supersolution property of $\underline{u}$ can be proved similarly.
As usual, we start from the partial dynamic programming principle: for $0\le
t<t+h\le T$,
\begin{eqnarray}  \label{DPP}
\overline{u}(t,x) &\le&\mathop{\overline{\rm lim}}_{n\rightarrow\infty}%
\mathbb{E}^G\Big[\int_t^{t+h} \delta_n(r) \alpha:d\langle B\rangle_r +
\overline{u}(t+h, x+ B_{t+h}^t)\Big],
\end{eqnarray}
Indeed, by the time homogeneity of the problem, we have
\begin{eqnarray*}
\overline{u}^n(t,x) &=& \mathbb{E}^G\Big[\int_t^{t+h} \delta_n(r)
\alpha:d\langle B\rangle_r + \mathbb{E}^G_{t+h}\big[ \int_{t+h}^T
\delta_n(r) \alpha: d\langle B\rangle_r + \varphi(x+B_T^{t})\big]\Big] \\
&=&\mathbb{E}^G\Big[\int_t^{t+h} \delta_n(r) \alpha: d\langle B\rangle_r +
\overline{u}^n(t+h,x+B_{t+h}^t)\Big]
\end{eqnarray*}
Then
\begin{eqnarray*}
&&\overline{u}(t,x) -\mathop{\overline{\rm lim}}_{n\rightarrow\infty}\mathbb{%
E}^G\Big[\int_t^{t+h} \delta_n(r) \alpha: d\langle B\rangle_r + \overline{u}%
(t+h, x+ B_{t+h}^t)\Big] \\
&=& \mathop{\overline{\rm lim}}_{n\rightarrow\infty} \overline{u}^n(t,x) -%
\mathop{\overline{\rm lim}}_{n\rightarrow\infty}\mathbb{E}^G\Big[%
\int_t^{t+h} \delta_n(r) \alpha: d\langle B\rangle_r + \overline{u}(t+h, x+
B_{t+h}^t)\Big] \\
&\le& \mathop{\overline{\rm lim}}_{n\rightarrow\infty}\mathbb{E}^G\Big[ (%
\overline{u}^n- \overline{u})(t+h, x+ B_{t+h}^t)\Big].
\end{eqnarray*}
Following standard arguments it is obvious that $\overline u$ is uniformly
Lipschitz continuous in $x$. Moreover, $\mathop{\overline{\rm lim}}_{n\to
\infty} (\overline{u}^n- \overline{u})(t+h,x) =0$ for any $x\in\mathbb{R}$.
Then \textrm{(\ref{DPP})} follows directly from the simple Lemma \ref%
{lem-GCT} below.

We next derive the viscosity subsolution property from \textrm{(\ref{DPP})}.
Let $(t, x)\in [0,T)\times \mathbb{R}^{d}$ and $\varphi\in
C^{1,2}([t,T)\times \mathbb{R}^{d})$ such that $0= [\overline{u}-\varphi](t,
x)=\max_{(s,y)\in [t,T]\times \mathbb{R}^d} [\overline u-\varphi](s,y)$.
Denote $X_s:= x+ B_s^{t}$. For any $0<h\le T-t$, by \textrm{(\ref{DPP})}
and then applying It\^{o}'s formula we have
\begin{eqnarray*}
&&\varphi(t,x) = \overline u(t,x) \le \mathop{\overline{\rm lim}}%
_{n\rightarrow\infty}\mathbb{E}^G\Big[\int_{t}^{t+h} \delta_n(r)\alpha:
d\langle B\rangle_r + \overline{u}(t+h,X_{t+h})\Big] \\
&\le&\mathop{\overline{\rm lim}}_{n\rightarrow\infty}\mathbb{E}^G\Big[%
\int_{t}^{t+h} \delta_n(r)\alpha: d\langle B\rangle_r + \varphi(t+h, X_{t+h})%
\Big] \\
&=&\mathop{\overline{\rm lim}}_{n\rightarrow\infty}\mathbb{E}^G\Big[%
\int_{t}^{t+h} \delta_n(r) \alpha:d\langle B\rangle_r+ \varphi(t,x)\\
&&  +
\int_{t}^{t+h} [\partial_t \varphi(r, X_r) dr + {\frac{1}{2}} \partial_{xx}
\varphi (r, X_r): d\langle B\rangle_r]\Big] \\
&\le& \mathop{\overline{\rm lim}}_{n\rightarrow\infty}\mathbb{E}^G\Big[%
\int_{t}^{t+h} [\alpha\delta_n(r) + {\frac{1}{2}} \partial_{xx} \varphi
(t,x)]:d\langle B\rangle_r\Big] + \varphi(t,x) + \partial_t\varphi(t,x)h \\
&&+ \mathbb{E}^G\Big[\int_{t}^{t+h} [\partial_t \varphi(r,
X_r)-\partial_t\varphi(t,x)] dr + {\frac{1}{2}}\int_t^{t+h}[\partial_{xx}
\varphi (r, X_r) - \partial_{xx}\varphi(t,x)] :d\langle B\rangle_r]\Big] \\
&\le& G^\alpha(\partial_{xx}\varphi(t,x)) h + \varphi(t,x) +
\partial_t\varphi(t,x)h \\
&&+ \mathbb{E}^G\Big[\sup_{t\le r\le t+h} [|\partial_t \varphi(r,
X_r)-\partial_t\varphi(t,x)|+ {\frac{|\overline\sigma^2|}{2}}|\partial_{xx}
\varphi (r, X_r) - \partial_{xx}\varphi(t,x)| \Big] h,
\end{eqnarray*}
thanks to \textrm{(\ref{dg})}. By standard arguments $\overline u$ is
uniformly Lipschitz continuous in $x$, and note that viscosity property is a
local property. Then, without loss of generality we may assume $\partial_t
\varphi$ and $\partial_{xx}$ are bounded and uniformly continuous in $(t,x)$
with a modulus of continuity function $\rho$. Thus,
\begin{eqnarray*}
0 \le \partial_t\varphi(t,x) + G^\alpha(\partial_{xx}\varphi(t,x)) + C%
\mathbb{E}^G\Big[\rho\big(C[h+\sup_{t\le r\le t+h} |B_r^t|]\big)\Big].
\end{eqnarray*}
Send $h\to 0$ we can easily get
\begin{eqnarray*}
\partial_t\varphi(t,x) + G^\alpha(\partial_{xx}\varphi(t,x)) \ge 0.
\end{eqnarray*}
Clearly $\overline u(T,x) = \varphi$. Therefore, $\overline u$ is a
viscosity subsolution of PDE \textrm{(\ref{ua})}. \hfill \vrule width.25cm
height.25cm depth0cm\smallskip

\begin{lem}
\label{lem-GCT}
Assume $\f_n: \dbR^d\to\dbR$ are uniformly Lipschitz continuous functions, uniformly in $n$, and $\limsup_{n\to\infty} \f_n(x)\le 0$ for all $x$. Then $\limsup_{n\to\infty} \dbE^G[\f_n(B_t)]\le 0$ for any $t$.
\end{lem}
{\noindent \textit{Proof.\quad}} Let $L$ denote the uniform Lipschitz
constant of $\varphi_n$. For any $\varepsilon>0$ and $R>0$, there exist
finitely many $x_i$, $i=1,\cdots, M$ and a partition $\cup_{i=1}^M O_i = O_R(%
\mathbf{0}) := \{x\in\mathbb{R}^d: |x|\le R\}$ such that $|x-x_i|\le
\varepsilon$ for all $x\in O_i$. Denote $O_0 := \mathbb{R}^d\backslash O_R(%
\mathbf{0})$ and $x_0:=\mathbf{0}$. Then
\begin{eqnarray*}
\varphi_n(B_t) &=& \sum_{i=0}^M \varphi_n(B_t)\mathbf{1}_{O_i}(B_t) =
\sum_{i=0}^M \varphi_n(x_i)\mathbf{1}_{O_i}(B_t)+ \sum_{i=0}^M
[\varphi_n(B_t) - \varphi_n(x_i)] \mathbf{1}_{O_i}(B_t) \\
&\le& \sum_{i=0}^M \varphi_n^+(x_i)\mathbf{1}_{O_i}(B_t) + L|B_t|\mathbf{1}%
_{O_0}(B_t) + L\varepsilon \sum_{i=1}^M\mathbf{1}_{O_i}(B_t) \\
&\le& \sum_{i=0}^M \varphi_n^+(x_i)\mathbf{1}_{O_i}(B_t) + {\frac{L}{R}}%
|B_t|^2 + L\varepsilon.
\end{eqnarray*}
Thus, noting that our condition implies $\limsup_{n\to\infty} \f^+_n(x)= 0$,
\begin{eqnarray*}
&&\mathop{\overline{\rm lim}}_{n\rightarrow\infty}\mathbb{E}^G\Big[%
\varphi_n(B_t)\Big] \le \mathop{\overline{\rm lim}}_{n\rightarrow\infty}%
\mathbb{E}^G\Big[ \sum_{i=0}^M \varphi_n^+(x_i)\mathbf{1}_{O_i}(B_t)+ {%
\frac{L}{R}}|B_t|^2 + L\varepsilon\Big] \\
&\le& \sum_{i=0}^M \mathop{\overline{\rm lim}}_{n\rightarrow\infty}%
\varphi_n^+(x_i)\mathbb{E}^G[\mathbf{1}_{O_i}(B_t)] + {\frac{L}{R}}
\mathbb{E}^G[|B_t|^2] + L\varepsilon
= {\frac{L}{R}} \mathbb{E}^G[|B_t|^2] + L\varepsilon.
\end{eqnarray*}
Send $R\to\infty$ and $\varepsilon\to 0$, we prove the result. \hfill \vrule %
width.25cm height.25cm depth0cm\smallskip

\subsection{An intermediate norm for $\protect\eta\in \mathcal{M}^1_G$}

We now use $\delta_n(t)$ to introduce the following norm for a
process $\eta$.
\begin{thm}
\label{thm-etaG}
For any $\eta\in\cM^1_G$, the following limit exists:
\bea
\label{dbMG}
\|\eta\|_{\dbM_G} := \lim_{n\to\infty}\dbE^G\Big[\int_0^T\delta_n(t)\eta_t : d\langle B\rangle_t\Big].
\eea
\end{thm}

{\noindent \textit{Proof.\quad}} We first assume $\eta \in \mathcal{M}^0_G$.
By otherwise considering a finer partition of $[0,T]$, without loss of
generality we assume, for $0=t_0<\cdots<t_m=T$,
\begin{eqnarray}  \label{eta}
\eta = \sum_{i=0}^{m-1} \eta_{t_i} \mathbf{1}_{[t_i, t_{i+1})}, &\mbox{where}%
& \eta_{t_i} = \varphi_i(B_{t_1},\cdots, B_{t_i})
\end{eqnarray}
and $\varphi_i$ is uniformly Lipschitz continuous. Denote
\begin{eqnarray*}
\psi^n_i (B_{t_1},\cdots, B_{t_i}) := \mathbb{E}^G_{t_i}\Big[%
\int_{t_i}^T\delta_n(t)\eta_t : d\langle B\rangle_t\Big].
\end{eqnarray*}
We prove by backward induction that
\begin{eqnarray}  \label{etaGi}
\lim_{n}\psi^n_i= \psi_i
\end{eqnarray}
where, $\psi_m:=0$ and, for $i=m-1,\cdots,0$,
\begin{eqnarray}  \label{psii}
\psi_i(x_1,\cdots, x_i) := \mathcal{E}^{{\varphi_i(x_1,\cdots,x_i)}}_{t_i,
t_{i+1}}(\psi_{i+1}(x_1,\cdots,x_i, \cdot)) (x_i).
\end{eqnarray}
Indeed, when $i=m$, \textrm{(\ref{etaGi})} holds obviously. Assume \textrm{(%
\ref{etaGi})} holds for $i+1$. Then by \textrm{(\ref{df})} we have
\begin{eqnarray*}
&& \mathop{\overline{\rm lim}}_{n\to\infty}\psi^n_i (B_{t_1},\cdots,
B_{t_i}) - \psi_i(B_{t_1},\cdots, B_{t_i}) \\
& =& \mathop{\overline{\rm lim}}_{n\to\infty}\mathbb{E}^G_{t_i}\Big[%
\int_{t_i}^{t_{i+1}}\delta_n(t) \eta_{t_i} : d\langle B\rangle_t +
\psi^n_{i+1}(B_{t_1},\cdots, B_{t_{i+1}})\Big] \\
&& - \lim_{n\to\infty}\mathbb{E}^G_{t_i}\Big[\int_{t_i}^{t_{i+1}}\delta_n(t)%
\eta_{t_i} : d\langle B\rangle_t + \psi_{i+1}(B_{t_1},\cdots, B_{t_{i+1}})%
\Big]\Big| \\
&\le& \mathop{\overline{\rm lim}}_{n\to\infty}\mathbb{E}^G_{t_i}\Big[%
\psi^n_{i+1}(B_{t_1},\cdots, B_{t_{i+1}})-\psi_{i+1}(B_{t_1},\cdots,
B_{t_{i+1}})\Big].
\end{eqnarray*}
By induction assumption, $\lim_{n\to\infty} \psi^n_{i+1} = \psi_{i+1}$.
Moreover, one can easily check that $\psi^n_{i+1}$ is uniformly continuous
in $x_{i+1}$, uniformly in $n$. Then by Lemma \ref{lem-GCT} we obtain
\begin{eqnarray*}
\mathop{\overline{\rm lim}}_{n\to\infty}\psi^n_i (B_{t_1},\cdots, B_{t_i}) -
\psi_i(B_{t_1},\cdots, B_{t_i}) \le 0.
\end{eqnarray*}
Similarly, we can show that
\begin{eqnarray*}
\psi_i(B_{t_1},\cdots, B_{t_i}) - \mathop{\underline{\rm lim}}%
_{n\to\infty}\psi^n_i (B_{t_1},\cdots, B_{t_i}) \le 0.
\end{eqnarray*}
Thus \textrm{(\ref{etaGi})} holds for $i$. This completes the induction and
hence proves that the limit in \textrm{(\ref{dbMG})} for $\eta\in \mathcal{M}%
^0_G$.

\medskip

We now consider general $\eta\in \mathcal{M}^1_G$. Let $\eta^m\in \mathcal{M%
}^0_G$ such that $\lim_{m\to\infty}\|\eta^m-\eta\|_{\mathcal{M}^1_G}=0$. For
each $m$, by previous arguments we have
\begin{eqnarray*}
\lim_{n\to\infty}\mathbb{E}^G\Big[\int_0^T\delta_n(t)\eta^m_t : d\langle
B\rangle_t\Big]~\mbox{exists}.
\end{eqnarray*}
By \reff{trest}, one can easily check that
\begin{eqnarray*}
&&\Big|\mathbb{E}^G\Big[\int_0^T\delta_n(t)\eta^m_t : d\langle B\rangle_t%
\Big] - \mathbb{E}^G\Big[\int_0^T\delta_n(t)\eta_t : d\langle B\rangle_t%
\Big] \Big| \\
&\le&\mathbb{E}^G\Big[\Big|\int_0^T\delta_n(t)[\eta^m_t - \eta_t] : d\langle
B\rangle_t\Big|\Big]\le \mathbb{E}^G\Big[\int_0^T|\eta^m_t - \eta_t|
|\overline \sigma^2|dt\Big] \\
&=& |\overline\sigma^2|\|\eta^m-\eta\|_{\mathcal{M}^1_G}\to 0,\quad\mbox{as}%
~m\to\infty.
\end{eqnarray*}
This clearly leads to the existence of $\lim_{n\to\infty}\mathbb{E}^G\Big[%
\int_0^T\delta_n(t)\eta_t : d\langle B\rangle_t\Big]$. \hfill \vrule %
width.25cm height.25cm depth0cm\smallskip

We now collect some basic properties of $\|\cdot\|_{\mathbb{M}_G}$. The left
inequality of \textrm{(\ref{etaGest})} below is crucial for our purpose. We
remark that, the norm $\|\cdot\|_{\mathcal{M}^1_{G_\varepsilon}}$ was
introduced by Hu and Peng \cite{HP2} and a similar estimate was obtained by
Song \cite{Song-unique} by using different arguments.  Recall the $c_0$ defined by \reff{c0}.

\begin{thm}
\label{thm-etaGproperty}
 $\|\cd\|_{\dbM_G}$ defines a norm on $\cM^1_G$, and for any $0<\e\le c_0$, it holds that,
 \bea \label{etaGest}
  \e \|\eta\|_{\cM^1_{G_\e}} \le \|\eta\|_{\dbM_G} \le C_0\|\eta\|_{\cM^1_G}.
  \eea
\end{thm}

To prove the theorem, we introduce some additional notations. Recall \textrm{(\ref{Ge})} and set
\begin{eqnarray}  \label{cPe}
\mathcal{A}_\varepsilon := \Big\{\sigma\in\mathcal{A}:
\underline\sigma_\varepsilon^2 \le \sigma^2 \le
\overline\sigma_\varepsilon^2 \Big\},\quad \mathcal{P}_\varepsilon := \Big\{%
\mathbb{P}^\sigma: \sigma\in \mathcal{A}_\varepsilon\Big\},\quad \mathcal{P}%
_0 := \lim_{\varepsilon\to 0} \mathcal{P}_\varepsilon.
\end{eqnarray}
We remark that the following inclusions are strict:
\bea
\label{cP0}
\mathcal{P}_0\subset  \{\mathbb{P}^\sigma:
\sigma \in \mathcal{A}, \underline\sigma< \sigma< \overline\sigma\}\subset \cP,~\mbox{but $\cP_0\subset \cP$ is dense under the weak topology}.
\eea

 \ms

{\noindent \textit{Proof.\quad}} (i) We first prove the estimates \reff{etaGest}. Note that $\|\cdot\|_{\mathcal{M}%
^1_{G_\varepsilon}} \le \|\cdot\|_{\mathcal{M}^1_G}$. By using standard
approximation arguments, it suffices to prove the statements for $\eta\in
\mathcal{M}^0_G$. We now assume $\eta$ takes the form \textrm{(\ref{eta})}
and we shall use the notations in the proof of Theorem \ref{thm-etaG}. In
particular, by \textrm{(\ref{etaGi})} we have
\begin{eqnarray*}
\|\eta\|_{\mathbb{M}_G} = \psi_0.
\end{eqnarray*}
%
Define $\underline\psi^%
\varepsilon_i $ and $\overline\psi^\varepsilon_i$ by:
\begin{eqnarray*}
\underline\psi^\varepsilon_i(B_{t_1},\cdots, B_{t_i}) := \varepsilon \mathbb{%
E}^{G_\varepsilon}_{t_i} \Big[\int_{t_i}^T |\eta_t|dt\Big],\quad
\overline\psi^\varepsilon_i(B_{t_1},\cdots, B_{t_i}) := C_0 \mathbb{E}%
^{G}_{t_i} \Big[\int_{t_i}^T |\eta_t|dt\Big].
\end{eqnarray*}
Then $\underline\psi^\varepsilon_m = \overline\psi^\varepsilon_m=0$, and
\begin{eqnarray*}
\underline\psi^\varepsilon_i(x_1,\cdots,x_i) &=& \mathbb{E}^{G_\varepsilon} %
\Big[\underline\psi^\varepsilon_{i+1}(x_1,\cdots, x_i, x_i + B_{t_{i+1}}^{t_i})%
\Big] +\varepsilon|\varphi_i(x_1,\cdots,x_i)|(t_{i+1}-t_i); \\
\overline\psi^\varepsilon_i(x_1,\cdots,x_i) &=& \mathbb{E}^{G_\varepsilon} %
\Big[\overline\psi^\varepsilon_{i+1}(x_1,\cdots, x_i, x_i + B_{t_{i+1}}^{t_i})%
\Big]+ C_0|\varphi_i(x_1,\cdots,x_i)|(t_{i+1}-t_i).
\end{eqnarray*}
Applying Lemma \ref{lem-Gae} (ii) and recalling \textrm{(\ref{psii})}, by
induction one proves \textrm{(\ref{etaGest})} immediately.

(ii) We now prove that $\|\cdot\|_{\mathbb{M}_G}$ defines a norm. Let $\eta\in \cM^1_G$. First, by 
\textrm{(\ref{etaGest})} we have $\|\eta\|_{\mathbb{M}_G}\ge 0$ and equality holds when $\eta=0$, $\cP$-q.s. On the other hand, assume $\|\eta\|_{\dbM_G}=0$, then by the left inequality of 
\textrm{(\ref{etaGest})} again we see that $\eta =0$, $\cP_0$-q.s. Now for any $\dbP\in \cP$, by \reff{cP0} there exists $\dbP_n\in \cP_0$ such that $\dbP_n$ converges to $\dbP$ weakly. Since $\eta\in \cM^1_G$, then $|\eta|$ is $\cP$-q.s. continuous and it follows from \cite{DHP} Lemma 27 that
\beaa
\dbE^\dbP\Big[\int_0^T |\eta_t|dt\Big] =\lim_{n\to\infty} \dbE^{\dbP_n}\Big[\int_0^T |\eta_t|dt\Big] =0.
\eeaa
That is, $\eta =0$, $\dbP$-a.s. for all $\dbP\in \cP$. Therefore, $\|\eta\|_{\mathbb{M}_G}= 0$ if and only if $\eta=0$, $\cP$-q.s.

Next,
for any $\lambda\in \mathbb{R}$, noting that $G^{-\alpha} = G^{\alpha}$ by
Lemma \ref{lem-Ga} (ii), it follows from \textrm{(\ref{psii})} that
\begin{eqnarray*}
\|\lambda\eta\|_{\mathbb{M}_G} = \||\lambda|\eta\|_{\mathbb{M}_G} &=&
\lim_{n\to\infty} \mathbb{E}^G\Big[\int_0^T \delta_n(t)
|\lambda|\eta_t:d\langle B\rangle_t\Big] \\
&=& \lim_{n\to\infty}|\lambda| \mathbb{E}^G\Big[\int_0^T \delta_n(t) \eta_t
: d\langle B\rangle_t\Big] = |\lambda| \|\eta\|_{\mathbb{M}_G}.
\end{eqnarray*}
Finally, for any $\eta, \tilde\eta\in \mathcal{M}^0_G$, by the sublinearity
of $\mathbb{E}^G$, we have
\begin{eqnarray*}
\mathbb{E}^G\Big[\int_0^T \delta_n(t) [\eta_t+\tilde\eta_t] : d\langle
B\rangle_t\Big] \le \mathbb{E}^G\Big[\int_0^T \delta_n(t) \eta_t : d\langle
B\rangle_t\Big] + \mathbb{E}^G\Big[\int_0^T \delta_n(t) \tilde\eta_t :
d\langle B\rangle_t\Big].
\end{eqnarray*}
Send $n\to\infty$ we obtain the triangle inequality: $\|\eta + \tilde\eta\|_{%
\mathbb{M}_G}\le \|\eta\|_{\mathbb{M}_G} + \|\tilde\eta\|_{\mathbb{M}_G}$.
That is, $\|\cdot\|_{\mathbb{M}_G}$ defines a
norm on $\cM^1_G$. \hfill \vrule width.25cm height.25cm depth0cm\smallskip

\subsection{The new norm for $\protect\eta$}

One drawback of the above norm $\|\cdot\|_{\mathbb{M}_G}$ is that we have to
use different norms in the left and right sides of \textrm{(\ref{etaGest})}. Consequently, we are not able to prove the completeness of $\mathcal{M}^1_G$
under $\|\cdot\|_{\mathbb{M}_G}$. To be precise, given a Cauchy sequence $%
\eta^n\in \mathcal{M}^1_G$ under $\|\cdot\|_{\mathbb{M}_G}$, we are not able
to prove the existence of a process $\eta$ such that $\displaystyle%
\lim_{n\to\infty}\|\eta^n-\eta\|_{\mathbb{M}_G}= 0$. For this reason, we
shall modify $\|\cdot\|_{\mathbb{M}_G}$ slightly by using heavily
the estimate \textrm{(\ref{etaGest})}. Set $\varepsilon_k := \frac{1}{1+k}%
c_0 $, $k\ge 1$, and define
\begin{eqnarray}  \label{dbM*G}
\|\eta\|_{\mathbb{M}^*_G}:= \sum_{k=1}^\infty 2^{-k} \|\eta\|_{\mathbb{M}%
_{G_{\varepsilon_k}}},\quad \eta\in \mathcal{M}^1_{G}.
\end{eqnarray}
Then clearly $\|\cdot\|_{\mathbb{M}^*_G}$ defines a norm on $\mathcal{M}^1_G$%
, and we denote by $\mathbb{M}^*_G$ the closure of $\mathcal{M}^1_G$ under $%
\|\cdot\|_{\mathbb{M}^*_G}$. To understand the space $\mathbb{M}^*_G$, we
note that $\mathcal{M}^1_{G_\varepsilon}$ is decreasing as $\varepsilon\to 0$. Set
\begin{eqnarray}  \label{cM1G0}
\mathcal{M}^1_{G_0} := \lim_{\varepsilon\to 0} \mathcal{M}^1_{G_\varepsilon}
= \bigcap_{0<\varepsilon \le c_0} \mathcal{M}^1_{G_\varepsilon} .
\end{eqnarray}

\begin{rem}
\label{rem-M*G}
{\rm (i) As mentioned earlier, elements in $\cM^1_G$ (resp. $\cM^1_{G_\e}$)  are considered identical if they are equal $\cP$-q.s. (resp. $\cP_\e$-q.s.).  Similarly elements in $\dbM^*_G$ are considered identical if they are equal $\cP_0$-q.s.

(ii) Obviously $\cM^1_{G_\e}\downarrow \cM^1_{G_0} $  as $\e\downarrow 0$. Thus the space $\cM^1_{G_0}$ is independent of $c_0$.

(iii) By \reff{etaGest}, it is obvious that
\bea
\label{dbM*Ginclusion}
\cM^1_G \subset \dbM^*_G \subset \cM^1_{G_0}.
\eea
Moreover, the above inclusions are strict. Indeed, consider the case $d=1$ for simplicity. One may easily see that $\eta_t := \1_{\{\la B\ra_t = \underline \si^2\}}$ is in $\dbM^*_G \backslash \cM^1_G $,  and $\eta_t := \sum_{n=1}^\infty  2^n \f_n\big({\la B\ra_t - \underline\si^2 t \over(\overline\si^2 - \underline\si^2) t}\big)$ is in  $\cM^1_{G_0}\backslash  \dbM^*_G$, where $\f_n$ is the linear interpolation such that $\f_n(\g) = 0$ when $\g \le {1\over n+1}$ or $\g \ge {1\over n}$, and $\f_n(\g) = 1$ when $\g = {1\over 2}[{1\over n}+{1\over n+1}]$ .
\qed}
\end{rem}
 We now
have
\begin{thm}
\label{thm-M*G}
The space $\dbM^*_G$ is complete under the norm $\|\cd\|_{\dbM^*_G}$. 
\end{thm}
{\noindent \textit{Proof.\quad}} First, it is clear that $\|\cd\|_{\dbM^*_G}$ is a seminorm on  $\dbM^*_G$.  Now assume $\eta\in \dbM^*_G$ such that $\|\eta\|_{\dbM^*_G} =0$. By  \reff{etaGest},  $\|\eta\|_{\cM^1_{G_\e}} =0$ for all $\e\le c_0$. Then $\eta =0$, $\cP_\e$-q.s. for all $0<\e\le c_0$ and thus $\eta =0$, $\cP_0$-q.s. That is, $\|\cd\|_{\dbM^*_G}$ is a norm on $\dbM^*_G$ (again, in the $\cP_0$-q.s. sense).

It remains to prove the completeness of the space. Let $\eta^n\in \dbM^*_{G}$ be a Cauchy sequence under $\|\cd\|_{\dbM^*_G}$.
For any $0<\varepsilon \le c_0$, there exists $k$ large enough such
that $\varepsilon_k< \varepsilon$. By the left inequality of \textrm{(\ref%
{etaGest})} we see that
\begin{eqnarray*}
\|\eta^n-\eta^m\|_{\mathcal{M}^1_{G_\varepsilon}} \le C_{\varepsilon,
\varepsilon_k} \|\eta^n-\eta^m\|_{\mathbb{M}_{G_{\varepsilon_k}}} \le 2^k
C_{\varepsilon, \varepsilon_k} \|\eta^n-\eta^m\|_{\mathbb{M}^*_{G}} \to
0,\quad\mbox{as}~n,m\to\infty.
\end{eqnarray*}
Since $(\mathcal{M}^1_{G_\varepsilon}, \|\cdot\|_{\mathcal{M}%
^1_{G_\varepsilon}})$ is complete, there exists unique (in $\mathcal{P}%
_\varepsilon$-q.s. sense) $\eta^{(\varepsilon)}\in \mathcal{M}%
^1_{G_{\varepsilon}}$ such that $\lim_{n\to\infty}\|\eta^n-\eta^{(%
\varepsilon)}\|_{\mathcal{M}^1_{G_\varepsilon}} = 0$. By the uniqueness,
clearly $\eta^{(\varepsilon)} = \eta^{(\tilde\varepsilon)}$, $\mathcal{P}%
_\varepsilon$-q.s. for any $0<\tilde\varepsilon<\varepsilon \le c_0$. Thus
there exists $\eta \in \mathcal{M}^1_{G_0}$ such that $\eta^{(\varepsilon)}
= \eta$, $\mathcal{P}_\varepsilon$-q.s. for all $0<\varepsilon \le c_0$.

We now show that
\bea
\label{etanconv}
\lim_{n\to\infty} \|\eta^n - \eta\|_{\dbM^*_G} = 0.
\eea
Indeed, for any $\delta>0$, there
exists $N_\delta$ such that
\begin{eqnarray*}
\|\eta^n-\eta^m\|_{\mathbb{M}^*_{G}} \le \delta,\quad \mbox{for all}~n, m\ge
N_\delta.
\end{eqnarray*}
Note that, by the right inequality of \textrm{(\ref{etaGest})},
\begin{eqnarray*}
\|\eta^m - \eta\|_{\mathbb{M}_{G_{\varepsilon}}} \le C_0\|\eta^m -
\eta^{(\varepsilon)}\|_{\mathcal{M}^1_{G_{\varepsilon}}} \to 0,\quad\mbox{as}%
~m\to\infty.
\end{eqnarray*}
Then for any $n\ge N_\delta$ and $K\ge 1$,
\begin{eqnarray*}
\sum_{k=1}^K 2^{-k} \|\eta^n-\eta\|_{\mathbb{M}_{G_{\varepsilon_k}}} =
\lim_{m\to\infty}\sum_{k=0}^K 2^{-k} \|\eta^n-\eta^m\|_{\mathbb{M}%
_{G_{\varepsilon_k}}} \le \mathop{\underline{\rm lim}}_{m\to\infty}\|\eta^n-%
\eta^m\|_{\mathbb{M}^*_{G}} \le \delta.
\end{eqnarray*}
Send $K\to\infty$ we obtain $\|\eta^n - \eta\|_{\mathbb{M}^*_G} \le \delta$
for all $n\ge N_\delta$. This proves \textrm{(\ref{etanconv})}, and hence the theorem.
%
\hfill \vrule %
width.25cm height.25cm depth0cm\smallskip

\section{The $G$-martingale representation theorem}

\label{sect-mrt}\setcounter{equation}{0} We first note that, assuming $%
(Y^i, Z^i, \eta^i)$, $i=1,2$, satisfy \textrm{(\ref{Gmg})}, then
\begin{eqnarray}
\label{DY}
d(Y^1_t - Y^2_t) = (Z^1_t-Z^2_t) \cdot dB_t - [G(\eta^1_t) - G(\eta^2_t)]dt + {\frac{1}{2}} [\eta^1_t - \eta^2_t]
: d\langle B\rangle_t.
\end{eqnarray}
By Lemma \ref{lem-detadt} we have
\begin{eqnarray}  \label{etaY}
\|\eta^1-\eta^2\|_{\mathbb{M}_{G_\varepsilon}} =2 \lim_{n\to\infty} \mathbb{E%
}^{G_\varepsilon}\Big[\int_0^T \delta_n(t) d(Y^1_t-Y^2_t)\Big],\quad %
\mbox{for all}~0<\varepsilon\le c_0.
\end{eqnarray}
In light of \textrm{(\ref{dbM*G})}, for any $p>1$ we define:
\begin{eqnarray}  \label{rhoG}
d_{G,p}(\xi_1,\xi_2):= \|Y^1-Y^2\|_{\mathbb{D}^p_G} + \sum_{k=1}^\infty
2^{-k} \lim_{n\to\infty} \mathbb{E}^{G_{\varepsilon_k}}\Big[\int_0^T
\delta_n(t) d(Y^1_t-Y^2_t)\Big], \\
\mbox{where}~\xi_i\in\mathcal{L}_{ip}~~\mbox{and}~~Y^i_t := \mathbb{E}%
^G_t[\xi_i] ,i=1,2.\qquad\qquad\qquad\qquad\qquad  \notag
\end{eqnarray}
Then clearly $d_{G,p}$ is a metric on $\mathcal{L}_{ip}$, and we let $\dbL^{*p}_G\subset \cL^p_G$ denote the closure of $\cL_{ip}$ under $d_{G,p}$. We remark that
\begin{eqnarray*}
\|\xi_1-\xi_2\|_{\mathcal{L}^p_G}\le \|Y^1-Y^2\|_{\mathbb{D}^p_G}\le
\|\xi_1-\xi_2\|_{\mathbb{L}^p_G}.
\end{eqnarray*}

\begin{rem}
\label{rem-mrt}
{\rm We remark that  we allow the metric $d_{G,p}(\xi_1,\xi_2)$ to depend on $Y^i$,  but not on $Z^i$ or $\eta^i$ explicitly. The component $Y$ has a representation, namely as the conditional $G$-expectation of $\xi$, but in general we do not have a desirable representation for $Z$ or $\eta$.  Thus it is relatively easier to check conditions imposed on $Y$ than those on $Z$ or $\eta$. See also \cite{PZ} for similar idea.
\qed}\end{rem}

Given $(Z, \eta)$ and $y$, let $Y^{y,Z,\eta}$ denote the $G$-martingale defined by \reff{Gmg} with initial value $Y_0=y$.
We first have
\begin{lem}
\label{lem-mrt}
For any $p>1$, $y\in \dbR$, and $(Z, \eta)\in \cH^p_G\times \mathcal{M}^p_G$, we have  $Y^{y,Z,\eta}_T\in \dbL^{*p}_G$. Moreover, for any such $(y_i, Z^i, \eta^i)$, $i=1, 2$, we have
\begin{eqnarray}
\label{rhoGest} d_{G,p}(Y^{y_1,Z^1,\eta^1}_T, Y^{y_2,Z^2,\eta^2}_T)
\le C_p\Big[|y_1-y_2| + \|Z^1-Z^2\|_{\cH^p_G} +
|\eta^1-\eta^2\|_{\cM^p_G}\Big].
\end{eqnarray}
\end{lem}

{\noindent \textit{Proof.\quad}} We first prove the a priori estimate \reff{rhoGest}. Denote $Y^i := Y^{y_i,Z^i,\eta^i}$, $i=1,2$. By \reff{DY}, it is obvious that
  \bea
  \label{rhoGest1}
  \|Y^1-Y^2\|_{\dbD^p_G} \le C_p\Big[|y_1-y_2| + \|Z^1-Z^2\|_{\cH^p_G} + \|\eta^1-\eta^2\|_{\cM^p_G}\Big].
  \eea
  Moreover, by \reff{etaY} and the right inequality of \reff{etaGest}, we have
  \beaa
 && \lim_{n\to\infty} \mathbb{E}^{G_{\varepsilon_k}}\Big[\int_0^T
\delta_n(t) d(Y^1_t-Y^2_t)\Big] = {1\over 2}\|\eta^1-\eta^2\|_{\mathbb{M}_{G_{\varepsilon_k}}}\\
 &\le& C\|\eta^1-\eta^2\|_{\cM^1_{G_{\varepsilon_k}}}\le C\|\eta^1-\eta^2\|_{\cM^1_G} \le C\|\eta^1-\eta^2\|_{\cM^p_G}.
\eeaa
Then,
\beaa
\sum_{k=1}^\infty
2^{-k} \lim_{n\to\infty} \mathbb{E}^{G_{\varepsilon_k}}\Big[\int_0^T
\delta_n(t) d(Y^1_t-Y^2_t)\Big] \le C\sum_{k=1}^\infty
2^{-k}\|\eta^1-\eta^2\|_{\cM^p_G} = C\|\eta^1-\eta^2\|_{\cM^p_G}.
\eeaa
 This, together with \reff{rhoGest1}, implies \reff{rhoGest}.

 We now show that $Y_T := Y^{y,Z,\eta}_T\in \dbL^{*p}_G$ in two steps.

 {\it Step 1.} Assume $\eta = 0$. By \reff{rhoGest} and the definition of $\cH^p_G$, we may assume without loss of generality that $Z = \sum_{i=0}^{n-1} Z_{t_i}\1_{[t_i, t_{i+1})}\in \cH^0_G$. Then
\beaa
Y_T = Y_0 + \sum_{i=0}^{n-1} Z_{t_i} B_{t_{i+1}}^{t_i} \in \cL_{ip}\subset \dbL^{*p}_G.
\eeaa

{\it Step 2.} For the general case, by \reff{rhoGest} and the
definition of $\cM^p_G$, we may assume without loss of generality
that $\eta = \sum_{i=0}^{n-1} \eta_{t_i}\1_{[t_i, t_{i+1})}\in
\cM^0_G$. Then \beaa Y_T = Y_0 + \int_0^T Z_t \cdot dB_t -
\sum_{i=0}^{n-1}\Big[ G(\eta_{t_i}) [t_{i+1}-t_i] - {1\over
2}\eta_{t_i}: [\la B\ra_{t_{i+1}}-\la B\ra_{t_i}]\Big]. \eeaa For
each $i$, applying It\^{o}'s formula we have \beaa d\Big(
B_t^{t_i}(B_t^{t_i})^T\Big) = 2 B_t^{t_i} d(B_t^{t_i})^T + d\la
B^{t_i}\ra_t = 2 B_t^{t_i} dB^T_t + d\la B\ra_t,~~ t\in [t_i,
t_{i+1}]. \eeaa Then \beaa \eta_{t_i}: [\la B\ra_{t_{i+1}}-\la
B\ra_{t_i}]\Big] = \eta_{t_i}:
[B_{t_{i+1}}^{t_i}(B_{t_{i+1}}^{t_i})^T] -
2\int_{t_i}^{t_{i+1}}(\eta_{t_i} B_t^{t_i}) \cdot dB_t. \eeaa Thus
\bea \label{YT} Y_T = Y_0 + \int_0^T \tilde Z_t \cdot dB_t -
\sum_{i=0}^{n-1}\Big[ G(\eta_{t_i}) [t_{i+1}-t_i] - {1\over
2}\eta_{t_i}: [B_{t_{i+1}}^{t_i}(B_{t_{i+1}}^{t_i})^T] \Big], \eea
where \beaa \tilde Z_t := Z_t -\sum_{i=0}^{n-1} \eta_{t_i} B^{t_i}_t
\1_{[t_i, t_{i+1})}(t). \eeaa One can easily check that $\tilde Z
\in \cH^p_G$. Then by Step 1,  $\int_0^T \tilde Z_t \cdot dB_t\in
\dbL^{*p}_G$. Moreover, it is obvious that
$\dis\sum_{i=0}^{n-1}\Big[ G(\eta_{t_i}) [t_{i+1}-t_i] - {1\over
2}\eta_{t_i}: [B_{t_{i+1}}^{t_i}(B_{t_{i+1}}^{t_i})^T] \Big]\in
\cL_{ip}$. Then it follows from \reff{YT} that $Y_T \in
\dbL^{*p}_G$. \qed

 Our main result of the paper is the following representation theorem, which is in the opposite direction of Lemma \ref{lem-mrt}.

\begin{thm}
\label{thm-mrt} Let $p>1$.

(i) For any $\xi\in \dbL^{*p}_G$ and denoting $Y_t := \dbE^G_t[\xi]$, there exist unique $Z\in \cH^p_G$ and $\eta\in \dbM^*_G$ such that \reff{Gmg} holds $\cP_0$-q.s. Moreover, there exists a constant $C_p>0$ such that
\bea
\label{mrtest1}
\|Y\|_{\dbD^p_G}+\|Z\|_{\cH^p_G} + \|\eta\|_{\dbM^*_G} \le C_pd_{G,p}(\xi,0).
\eea

(ii) For any $\xi_1, \xi_2\in \dbL^{*p}_G$, let $(Y^i, Z^i, \eta^i)$ denote the corresponding terms. Then
\bea
\label{mrtest2}
\|Y^1-Y^2\|_{\dbD^p_G} + \|\eta^1-\eta^2\|_{\dbM^*_G} \le C_pd_{G,p}(\xi_1,\xi_2),\q  \|Z^1-Z^2\|_{\cH^p_G} \le C_p\Big(d_{G,p}(\xi_1,\xi_2)\Big)^{1\over 2}.
\eea
\end{thm}

\begin{rem}
\label{rem-mrt2}
{\rm We can only prove the representation \reff{Gmg} in $\cP_0$-q.s.  sense. This is mainly because we are not able to prove the equivalence of $\|\cd\|_{\dbM_G}$ and $\|\cd\|_{\cM^1_G}$ in Theorem \ref{thm-etaGproperty}. See also Remark \ref{rem-M*G} (iii). It is still an open problem to establish the representation  \reff{Gmg}  $\cP$-q.s.  and we shall leave it for future research.
\qed}
\end{rem}
{\noindent \textit{Proof.\quad}} We proceed in two steps.

\textit{Step 1.} We first prove a priori estimates \textrm{(\ref{mrtest1})}
and \textrm{(\ref{mrtest2})} by assuming $(Y, Z, \eta)$ and $(Y^i, Z^i,
\eta^i)$, $i=1,2$, are in $\mathbb{D}^p_G\times \mathcal{H}^p_G\times
\mathbb{M}^*_G$ and satisfy \textrm{(\ref{Gmg})} $\mathcal{P}_0$-q.s.
Indeed, by \textrm{(\ref{etaY})} and \textrm{(\ref{rhoG})} it is clear that
\begin{eqnarray*}
\|Y\|_{\mathbb{D}^p_G} + \|\eta\|_{\mathbb{M}^*_G} \le C_pd_{G,p}(\xi,0),
\quad \|Y^1-Y^2\|_{\mathbb{D}^p_G} + \|\eta^1-\eta^2\|_{\mathbb{M}^*_G} \le
Cd_{G,p}(\xi_1,\xi_2).
\end{eqnarray*}
Moreover, combining the arguments in \cite{STZ-G} and \cite{HP1}, or
following the arguments in \cite{Song-G}, one can easily prove
\begin{eqnarray*}
\|Z\|_{\mathcal{H}^p_G} \le C_p\|Y\|_{\mathbb{D}^p_G},\quad \|Z^1-Z^2\|_{%
\mathcal{H}^p_G} \le C_p\Big( \|Y^1-Y^2\|_{\mathbb{D}^p_G} \Big)^{\frac{1}{2}%
}.
\end{eqnarray*}
Then \textrm{(\ref{mrtest1})} and \textrm{(\ref{mrtest2})} hold.

\textit{Step 2.} We next prove the existence of $(Z,\eta)$. For any $\xi\in \mathbb{L}^{*p}_G$, by definition there
exist $\xi_n\in \mathcal{L}_{ip}$ such that $\dis\lim_{n\to\infty}\rho_G^p(%
\xi_n,\xi) = 0$. Let $(Y^n, Z^n, \eta^n)$ be corresponding to $\xi_n$. As $%
n,m\to\infty$, by \textrm{(\ref{mrtest2})} we have
\begin{eqnarray*}
\|Y^n-Y^m\|_{\mathbb{D}^p_G} + \|\eta^n-\eta^m\|_{\mathbb{M}^*_G}+
\|Z^n-Z^m\|_{\mathcal{H}^p_G} \le C_p\Big[d_{G,p}(\xi_n,\xi_m)+\big(%
d_{G,p}(\xi_n,\xi_m)\big)^{\frac{1}{2}}\Big] \to 0.
\end{eqnarray*}
Then there exist $(Y, Z, \eta)\in \mathbb{D}^p_G\times \mathcal{H}^p_G\times
\mathbb{M}^*_G$ such that
\begin{eqnarray*}
\|Y^n-Y\|_{\mathbb{D}^p_G} + \|\eta^n-\eta\|_{\mathbb{M}^*_G}+ \|Z^n-Z\|_{%
\mathcal{H}^p_G} \to 0,\quad\mbox{as}~n\to\infty.
\end{eqnarray*}
Moreover, for any $0<\varepsilon\le c_0$, choose $k$ large enough so that $%
\varepsilon_k < \varepsilon$. Then
\begin{eqnarray*}
\|\eta^n- \eta\|_{\mathcal{M}^1_{G_\varepsilon}} \le C_{\varepsilon,
\varepsilon_k} \|\eta^n- \eta\|_{\mathbb{M}^1_{G_{\varepsilon_k}}} \le
2^kC_{\varepsilon, \varepsilon_k} \|\eta^n- \eta\|_{\mathbb{M}^*_G} \to
0,\quad\mbox{as}~n\to\infty.
\end{eqnarray*}
Thus
\begin{eqnarray*}
\int_0^t G(\eta^n_s)ds \to \int_0^t G(\eta_s)ds,\quad \mathcal{P}_\varepsilon%
\mbox{-q.s.}
\end{eqnarray*}
Since $(Y^n, Z^n, \eta^n)$ satisfy \textrm{(\ref{Gmg})} $\mathcal{P}%
_\varepsilon$-q.s., then it is clear that $(Y, Z, \eta)$ also satisfy
\textrm{(\ref{Gmg})} $\mathcal{P}_\varepsilon$-q.s. By the arbitrariness of $%
\varepsilon$ we see that $(Y, Z, \eta)$ also satisfy \textrm{(\ref{Gmg})} $%
\mathcal{P}_0$-q.s.

Finally, the uniqueness of $(Z, \eta)\in \mathcal{H}^p_G\times \mathbb{M}%
^*_G $ follows from \textrm{(\ref{mrtest2})}. \hfill \vrule width.25cm
height.25cm depth0cm\smallskip

We conclude this paper by providing a nontrivial example of $\xi$ which has the representation, but is not in $\mathcal{L}_{ip}$.

\begin{eg}
\label{eg-M*} Let $d=1$ and $B^*_t := \sup_{0\le s\le t}B_s$. Then $
B^*_T \in \dbL^{*p}_G\backslash \cL_{ip}$ for any $p>1$.
\end{eg}

{\noindent \textit{Proof.\quad}} It is clear that $B^*_T\notin \mathcal{L}%
_{ip}$. We prove $B^*_T \in \mathbb{L}^{*p}_G$ in several steps.

\textit{Step 1.} Assume $\xi: \O\to \dbR$ is uniformly Lipschitz continuous and convex in $\o$. We show that $\dbE^G[\xi] = \dbE^{\overline\dbP}[\xi]$, where $\overline \dbP := \dbP^{\overline\si}$.

Indeed, for any $n$, denote $t^n_i := {iT\over n}$, $i=0,\cds, n$,  $x_0 := 0$, and define
\beaa
g_n(x_1,\cds, x_n) &:=& \xi\Big(\sum_{i=1}^n {1\over t^n_i-t^n_{i-1}} [x_{i-1}(t^n_i - t)  + x_i (t-t^n_{i-1})]\1_{(t^n_{i-1}, t^n_i]}(t)\Big);\\
 \xi_n &:=& g_n(B_{t^n_1},\cds, B_{t^n_n}).
\eeaa
Since $\xi$ is convex, clearly $g_n$ is convex. Then
$
\dbE^G[\xi_n] = \dbE^{\overline \dbP}[\xi_n].
$
Since $\xi$ is uniformly Lipschitz continuous, then
\beaa
|\xi_n - \xi| \le C\max_{1\le i\le n} \sup_{t^n_{i-1}\le t\le t_i^n}|B_t-B_{t^n_i}|.
\eeaa
This implies that $\dbE^G[|\xi_n - \xi|\to 0$ and $\dbE^{\overline \dbP}[|\xi_n - \xi|]\to 0$ as $n\to \infty$, and therefore,  $\dbE^G[\xi] = \dbE^{\overline\dbP}[\xi]$.

{\it Step 2.} For simplicity, we assume $\overline{\sigma}=1$, and thus $\overline\dbP=\dbP_0$. Note that $\xi := B^*_T$ is uniformly Lipschitz continuous and convex in $\o$. Then by adapting Step 1 to conditional $G$-expectations we have
\beaa
Y_t := \dbE^G_t[\xi] = \dbE^{\dbP_0}_t [B^*_T] = u(t, B_t, B^*_t),
\eeaa
where, for $x\le y$,
\beaa
u(t, x, y) := \dbE^{\dbP_0}\Big[ y \vee [x + \sup_{t\le s\le T} B^t_s]\Big] = \dbE^{\dbP_0}\Big[ y \vee [x + B^*_{T-t}]\Big].
\eeaa
Note that, under $\dbP_0$, $B^*_{T-t}$ has the same distribution as $|B_{T-t}|$. Then
\beaa
u(t, x, y)&=&\sqrt{\frac{2}{\pi}}\int_0^\infty y\vee(x+\sqrt{T-t}z) e^{-\frac{z^2}{2}}dz \\
&=&\sqrt{\frac{2}{\pi}}\int_0^{\frac{y-x}{\sqrt{T-t}}} y e^{-\frac{z^2}{2}}dz +%
\sqrt{\frac{2}{\pi}}\int_{\frac{y-x}{\sqrt{T-t}}}^\infty(x+\sqrt{T-t}z) e^{-%
\frac{z^2}{2}}dz.
\eeaa
For $t\in[0,T)$ and $x<y$, we have
\bea
\label{pau}
&&\partial_t u(t,x, y) = -\frac{1}{\sqrt{2\pi (T-t)}}e^{-\frac{(y-x)^2}{2(T-t)}}; \q \partial_yu(t,x,y) = \sqrt{\frac{2}{\pi}}\int_0^{\frac{y-x}{\sqrt{T-t}}}e^{-\frac{z^2}{2}}dz;\nonumber\\
&&\partial_x u(t,x,y) = \sqrt{\frac{2}{\pi}}\int_{\frac{y-x}{\sqrt{T-t}}}^\infty
e^{-\frac{z^2}{2}}dz;\q
\partial_{xx} u(t,x,y) = \sqrt{\frac{2}{\pi (T-t)}}e^{-\frac{(y-x)^2}{2(T-t)}}>0.
\eea
Then
\beaa
\partial_t u+\frac{1}{2}G(\partial_{xx}u) = \partial_t u+\frac{1}{2}\partial_{xx}u =0,
&\mbox{and}&
\partial_yu (t,y,y)=0.
\eeaa
Note that $dB^*_t$ has support on $\{t: B^*_t = B_t\}$. Then by It\^o's formula we have
\begin{eqnarray*}
&&d Y_t= du(t,B_t, B^*_t) \\
&=& \pa_t u(t,B_t, B^*_t) dt + \pa_x u(t,B_t, B^*_t) dB_t + \pa_y u(t,B_t, B^*_t) dB^*_t + {1\over 2} \pa_{xx} u(t,B_t, B^*_t) d\la B\ra_t\\
&=& \pa_x u(t,B_t, B^*_t) dB_t - G(\pa_{xx} u(t,B_t, B^*_t)) dt + {1\over 2} \pa_{xx} u(t,B_t, B^*_t) d\la B\ra_t.
\end{eqnarray*}
Thus we obtain the representation with
\bea
\label{Zeta}
Z_t = \pa_x u(t,B_t, B^*_t),\q \eta_t = \pa_{xx} u(t,B_t, B^*_t).
\eea

\textit{Step 3.} By Lemma \ref{lem-mrt}, it remains to show that $(Z, \eta) \in \cH^p_G\times \cM^p_G$. For any $n$, denote
\beaa
Z^n_t := Z_t\1_{[0,T-{1\over n}]}, \q\eta^n_t := \eta_t\1_{[0, T-{1\over n}]}.
 \eeaa
 Note that, in the interval $[0, T-{1\over n}]$, $\pa_x u$ and $\pa_{xx} u$ are bounded and uniformly Lipschitz continuous in $(t,x,y)$, then clearly $(Z^n, \eta^n) \in \cH^p_G\times \cM^p_G$.
  Moreover, by \reff{pau} we have $|\pa_xu(t,x,y)|\le 1$ and $|\pa_{xx} u(t,x,y)|\le {C\over \sqrt{T-t}}$. Then, as $n\to\infty$,
 \beaa
 &&\dbE^G\Big[\Big(\int_0^T |Z_t - Z^n_t|^2 d\la B\ra_t\Big)^{p\over 2}\Big] =  \dbE^G\Big[\Big(\int_{T-{1\over n}}^T |Z_t|^2 d\la B\ra_t\Big)^{p\over 2}\Big] \\
 &\le& \dbE^G\Big[\Big(\la B\ra_T-\la B\ra_{T-{1\over n}}\Big)^{p\over 2}\Big] = {C_p\over n^{p\over 2}} \to 0;\\
 &&\dbE^G\Big[\Big(\int_0^T |\eta_t - \eta^n_t| dt\Big)^p\Big] =  \dbE^G\Big[\Big(\int_{T-{1\over n}}^T |\eta_t| dt\Big)^p\Big] \\
 &\le& C\dbE^G\Big[\Big(\int_{T-{1\over n}}^T {dt\over \sqrt{T-t}}\Big)^p\Big] = {C_p\over n^{p\over 2}} \to 0;
 \eeaa
This proves that $(Z, \eta) \in \cH^p_G\times \cM^p_G$ and completes the proof.
\qed

\providecommand{\bysame}{\leavevmode\hbox
to3em{\hrulefill}\thinspace} \providecommand{\MR}{\relax\ifhmode\unskip%
\space\fi MR }
\providecommand{\MRhref}[2]{  \href{http://www.ams.org/mathscinet-getitem?mr=#1}{#2}
} \providecommand{\href}[2]{#2}


\begin{thebibliography}{99}

\bibitem{CSTV}
Cheridito, P., Soner, H.M., Touzi, N., and Victoir, N. (2007)
{\it Second order BSDE's and fully nonlinear PDE's}, {\sl Communications
in Pure and Applied Mathematics}, 60, 1081-1110.

\bibitem{DHP} Denis, L., Hu, M. and Peng S. (2011) \textit{Function spaces
and capacity related to a sublinear expectation: application to $G$-Brownian
motion pathes}, \textsl{Potential Analysis}, 34 (2), 139-161.

\bibitem{DM} Denis, L. and Martini,  C. (2006) \textit{A Theorectical Framework
for the Pricing of Contingent Claims in the Presence of Model Uncertainty},
\textsl{Annals of Applied Probability}, 16 (2), 827-852.

\bibitem{DNS}
Dolinsky, Y., Nutz, M. and Soner, H. M. (2011) {\it Weak Approximation of G-Expectations},  preprint, arXiv:1103.0575.

\bibitem{ETZ}
Ekren, I., Touzi, N. and Zhang J. (2012) {\it Optimal Stopping under Nonlinear Expectation},  preprint, arXiv:1209.6601.

\bibitem{EJ}
Epstein, L. and Ji, S. (2011) {\it  Ambiguous Volatility, Possibility and Utility in Continuous Time},  preprint, arXiv:1103.1652.

\bibitem{FTW} Fahim, A., Touzi, N.  and Waxin, X. (2011) {\it A Probabilitic Numerical Method for Fully Nonlinear Parabolic PDEs}, {\sl Ann. Appl. Probab.}, 21(4), 1322-1364.

\bibitem{GZZ}
Guo, W., Zhang, J. and Zhuo, J. (2012) {\it A Monotone Scheme  for High Dimensional Fully Nonlinear  PDEs},  preprint, arXiv:1212.0466.

\bibitem{HP1} Hamadene, S. and Popier, A. (2008) \textit{$L^p$-Solutions for
Reflected Backward Stochastic Differential Equations},
preprint, arXiv:0807.1846.

\bibitem{HJPS}
Hu, M., Ji, S., Peng, S. and Song, Y. (2012) {\it Backward Stochastic Differential Equations Driven by G-Brownian Motion}, preprint, arXiv:1206.5889

\bibitem{HP2} Hu, Y. and Peng, S. (2010) \textit{Some Estimates for
Martingale Representation under $G$-Expectation}, preprint, arXiv:1004.1098.

\bibitem{MPZ} Matoussi, A., Possamai D. and  Zhou C. (2012) {\it Second Order Reflected Backward Stochastic Differential Equations}, preprint,
 arXiv:1201.0746.


\bibitem{Nutz}
 Nutz, M. (2010) {\it Random G-Expectations}, {\sl Annals of Applied Probability}, to appear.



\bibitem{NS}
Nutz, M. and Soner, H.M. (2012) {\it Superhedging and Dynamic Risk Measures under Volatility Uncertainty}, {\sl SIAM Journal on Control and Optimization}, to appear.

\bibitem{NV}
Nutz, M. and Van Handel, R. (2012) {\it Constructing Sublinear Expectations on Path Space}, preprint,  arXiv:1205.2415. 


\bibitem{NZ}
Nutz, M. and Zhang, J. (2012)  {\it Optimal Stopping under Adverse Nonlinear Expectation and Related Games}, preprint, arXiv:1212.2140.


\bibitem{Peng-g} Peng, S. (2004) \textit{Filtration consistent nonlinear
expectations and evaluations of contingent claims}, \textsl{Acta
Mathematicae Applicatae Sinica}, 20(2), 191-214.

\bibitem{Peng-G} Peng, S. (2006) \textit{$G$-Expectation, $G$-Brownian motion and related stochastic calculus of It\^{o}'s type},
\textsl{Stochastic Analysis and Applications, The Abel Symposium
2005, Abel Symposia}, Edit. Benth et. al., 541-567, Springer-Verlag.



\bibitem{Peng-book} Peng, S. (2010) \emph{Nonlinear Expectations and
Stochastic Calculus under Uncertainty}, preprint, arXiv:1002.4546.

\bibitem{Peng-ICM} Peng, S. (2010) \textit{Backward Stochastic Differential
Equation, Nonlinear Expectation and Their Applications}, \textsl{Proceedings
of the International Congress of Mathematicians 2010}, 393-432.

\bibitem{PZ} Pham, T. and Zhang, J. (2011) \textit{Some Norm Estimates for
Semimartingales -- Under Linear and Nonlinear Expectations}, preprint, arXiv:1107.4020v1.

\bibitem{STZ-G} Soner, M., Touzi, N., and Zhang, J. (2011) \textit{%
Martingale Representation Theorem under $G$-expectation}, \textsl{Stochastic
Processes and their Applications}, 121, 265-287.

\bibitem{STZ-Duality} Soner, M., Touzi, N., and Zhang, J. (2012) \textit{Dual Formulation of Second Order Target Problems}, \textsl{Annals of Applied Probability}, to appear.


\bibitem{STZ-2BSDE} Soner, M., Touzi, N., and Zhang, J. (2012) \textit{The
wellposedness of second order backward SDEs}, \textsl{Probability Theory and
Related Fields}, 153,  149-190.

\bibitem{Song-G} Song, Y.(2011) \textit{Some properties on $G$-evaluation
and its applications to $G$-martingale decomposition. } \textsl{Science
China Mathematics}, 54 (2), 287-300.

\bibitem{Song-unique} Song, Y.(2012) \textit{Uniqueness of the
representation for $G$-martingales with finite variation.}
Electron. J. Probab. 17, 1-15.

\bibitem{Tan}
Tan, X. (2012) {\it Probabilistic Numerical Approximation for Stochastic Control Problems}, preprint.

\bibitem{XZ} Xu, J. and Zhang B. (2009), \textit{Martingale characterization
of $G$-Brownian motion}, \textsl{Stochastic Processes and their Applications}, 119 (1), 232-248.
\end{thebibliography}
\end{document}